# On the Morse index of higher-dimensional free boundary minimal catenoids.

3rd September 2017

Graham Smith[1], Ari Stern[2], Hung Tran[3], Detang Zhou[4]

**Abstract:** For all $n$, we define the $n$-dimensional critical catenoid $M_n$ to be the unique rotationally symmetric, free boundary minimal hypersurface of non-trivial topology embedded in the closed unit ball in $\mathbb{R}^{n+1}$. We show that the Morse index $\mathrm{MI}(n)$ of $M_n$ satisfies the following asymptotic estimate as $n$ tends to infinity.

$$\underset{n\rightarrow+\infty}{\mathrm{Lim}}\frac{\mathrm{Log}(\mathrm{MI}(n))}{\sqrt{n}\mathrm{Log}(\sqrt{n})}=1.$$

We also study the numerical problem, providing exact values for the Morse index for $n=2,\cdots,100$, together with qualitative studies of $\mathrm{MI}(n)$ and related geometric quantities for large values of $n$.

[1] Instituto de Matemática, UFRJ, Av. Athos da Silveira Ramos 149, Centro de Tecnologia - Bloco C, Cidade Universitária - Ilha do Fundão, Caixa Postal 68530, 21941-909, Rio de Janeiro, RJ - BRAZIL

[2] Department of Mathematics, Washington University, Campus Box 1146, One Brookings Drive, St. Louis, MO

[3] Department of Mathematics and Statistics, Texas Tech University, Lubbock TX 79413.

[4] Instituto de Matemática, UFF, Rua Professor Marcos Waldemar de Freitas Reis, Bloco H - Campus do Gragoatá, São Domingos, 24.210-201, Niterói, RJ - BRAZIL

## 1 - Introduction.

**1.1 - Free boundary minimal catenoids.** Let $B := B^{n+1}$ denote the closed unit ball in $(n+1)$-dimensional Euclidean space. Let $M^n$ be a compact hypersurface smoothly embedded in $B$ with smooth boundary in $\partial B$. We say that $M^n$ is *free boundary minimal* whenever it is a critical point of the volume functional amongst all such embedded hypersurfaces. It is well known that this implies that $M^n$ has vanishing mean curvature, that $M^n$ meets $\partial B$ orthogonally along $\partial M^n$, and that no interior point of $M^n$ meets $\partial B$. Free boundary minimal surfaces, which have been studied since the classical work of Gergonne, have received renewed attention with the recent work [10] and [11] of Fraser & Schoen, and free boundary minimal hypersurfaces form a natural higher-dimensional extension of this class (c.f. [1], [2], [19] and [20]).

In this paper, we study questions of stability for certain free boundary minimal hypersurfaces. More precisely, suppose that $M^n$ is orientable and choose a unit normal vector field $N : M^n \to S^n$ which is compatible with the orientation. We identify infinitesimal perturbations of $M^n$ with smooth functions over this hypersurface by identifying the field $\phi N$ with the function $\phi$. Given $\phi$, the second order variation of volume for an infinitesimal perturbation in the direction of $\phi N$ is then given by (c.f. [16])

$$Q^n(\phi,\phi) := \int_M \|\nabla\phi\|^2 - \mathrm{Tr}(A^2)\phi^2 \mathrm{dV}^M - \int_{\partial M} \phi^2 \mathrm{dV}^{\partial M}, \tag{1.1}$$

where $\nabla$ here denotes the gradient operator of $M^n$, $A$ its shape operator, $\mathrm{dV}^M$ its volume element and $\mathrm{dV}^{\partial M}$ the volume element of $\partial M^n$. The *stability spectrum* of $M^n$ is defined to be the set of all eigenvalues in $\mathbb{C}$ of this symmetric bilinear form, where a function $\phi : M^n \to \mathbb{R}$ is said to be an *eigenfunction* of $Q^n$ with *eigenvalue* $\lambda$ whenever

$$Q^n(\phi,\psi) = \lambda\langle\phi,\psi\rangle \tag{1.2}$$

for any other function $\psi$, where $\langle\cdot,\cdot\rangle$ here denotes the $L^2$ inner product of $M^n$. It follows by the classical theory of elliptic operators (c.f. [12] and [15]) that the stability spectrum is discrete and real, and all but finitely many eigenvalues are positive. The *Morse index* of $M^n$ is then defined to be the number of eigenvalues, counted with multiplicity, that are strictly negative. Geometrically, this corresponds to the degree of *instability* of the hypersurface, in the sense that it is equal to the number of perturbation directions which reduce the volume.

Consider now the case where $M^n$ is the *critical catenoid*, which is defined to be the unique free boundary minimal hypersurface of revolution in $\mathbb{B}^{n+1}$ which has non-trivial topology. In previous work (c.f. [9], [21] and [24]), it was shown that $M^2$ has Morse index equal to 4, thereby solving a problem left open by Fraser & Schoen in [10] and [11]. In the current paper, we study the higher-dimensional case. Here, it turns out to be very hard to obtain a closed formula for the Morse index in general, and for this reason we satisfy ourselves with an asymptotic formula as $n$ tends to infinity. We show

**Theorem I**

*For all $n$, let $MI(n)$ denote the Morse index of the $n$-dimensional critical catenoid. Then*

$$\mathop{Lim}_{n\rightarrow+\infty}\frac{Log(MI(n))}{\sqrt{n}Log(\sqrt{n})}=1. \tag{1.3}$$

Theorem I follows immediately from Theorem 3.6.5, below.

We find it enlightening to compare Theorem I with the corresponding result for minimal hypersurfaces embedded in $S^{n+1}$, the unit sphere in $(n+2)$-dimensional Euclidean space. Indeed, for all $(p,q)$ such that $(p+q)=n$, let $M^{p,q}$ be the minimal hypersurface given by

$$M^{p,q}:=\sqrt{\frac{p}{p+q}}S^p\times\sqrt{\frac{q}{p+q}}S^q. \tag{1.4}$$

Here, if $N$ denotes the unit normal vector field over $M^{p,q}$ and if $\phi:M^{p,q}\rightarrow\mathbb{R}$ is some smooth function, then the second order variation of volume for an infinitesimal perturbation in the direction of $\phi N$ is given by

$$Q^{p,q}(\phi,\phi)=\int_M\|\nabla\phi\|^2-(n+\mathrm{Tr}(A^2))\phi^2\mathrm{dV}^M, \tag{1.5}$$

where $\nabla$ now denotes the gradient operator of $M^{p,q}$, $A$ denotes its shape operator, and $\mathrm{dV}^M$ denotes its volume element. As before, the *Morse index* of $M^{p,q}$ is defined to be the number of strictly negative eigenvalues of this symmetric bilinear form and in Appendix A we show that this number is equal to $(n+3)$. This has the following geometric interpretation. If $X$ is the vector field of an infinitesimal translation of $\mathbb{R}^{n+2}$, then the inner product $\phi:=\langle X,N\rangle$ is an eigenfunction of $Q^{p,q}$ with eigenvalue $(-n)$. Likewise, if $X$ is the vector field of an infinitesimal dilatation of either the $\mathbb{R}^{p+1}$ or the $\mathbb{R}^{q+1}$ component of $\mathbb{R}^{n+2}$, then the inner product $\phi:=\langle X,N\rangle$ is also an eigenfunction with eigenvalue $(-2n)$. Finally, since the space spanned by all such functions has dimension $(n+3)$, we see that this accounts for all negative eigenvalues of $Q^{p,q}$.

We now return our attention to free boundary minimal hypersurfaces in the unit ball. It is readily shown (c.f. [11]) that $Q^n$ is negative definite over the linear space $U^n$ generated by functions of the form $\phi:=\langle X,N\rangle$ where $X$ is the vector field of an infinitesimal translation or of an infinitesimal dilatation of $\mathbb{R}^{n+1}$ - both of which are known as *Killing fields*. This space has dimension $(n+2)$ and, although its elements are no longer eigenfunctions, the technique of Rayleigh quotients shows that the Morse index of $Q^n$ is nonetheless bounded below by this number.

In geometric terms, the result proven in [9], [21] and [24] states that elements of $U^2$ account for all negative eigenvalues of $Q^2$. Numerical methods then show that this continues to hold for $n=3$ and $n=4$, but not for $n\geq5$ (c.f. Table 10). Theorem I thus tells us that, as the dimension grows, $M^n$ has an ever greater number of unstable directions that are unrelated to Killing fields of the ambient space.

**1.2 - Steklov eigenvalues.** Theorem I may be reformulated into equivalent results which serve to illustrate the main ideas of the proof. First, the *Jacobi operator* $\mathrm{J}^n : C^\infty(M^n) \rightarrow C^\infty(M^n)$ is defined by

$$\mathrm{J}^n\phi := -\Delta\phi - \mathrm{Tr}(A^2)\phi, \tag{1.6}$$

where $\Delta$ here denotes the Laplace-Beltrami operator of $M^n$. The subspace $V$ of $C^\infty(\partial M^n)$ is then defined by

$$V := \{\partial_\nu\phi \mid \phi \in \mathrm{Ker}(\mathrm{J}^n)\}^\perp, \tag{1.7}$$

where the orthogonal complement is taken with respect to the $L^2$ inner product. Standard results of elliptic theory (c.f. [12]) show that for all $\phi \in V$, there exists a unique $\hat{\phi} \in \mathrm{Ker}(\mathrm{J}^n)^\perp$, known as the *Dirichlet extension* of $\phi$, such that

$$\hat{\phi} \circ e = \phi, \tag{1.8}$$

where $e : \partial M^n \rightarrow M^n$ here denotes the canonical embedding. The *Dirichlet to Neumann* operator $\mathrm{DN} : V \rightarrow V$ is then defined by

$$\mathrm{DN}\phi := \pi(\partial_\nu\hat{\phi}), \tag{1.9}$$

where $\pi : C^\infty(\partial M) \rightarrow V$ here denotes the orthogonal projection, again with respect to the $L^2$ inner product. The eigenvalues of this self-adjoint operator are known as the *Steklov eigenvalues* of $M^n$, and the *Steklov index* of $M^n$, which we denote by $\mathrm{SI}(n)$, is defined to be the number of its Steklov eigenvalues, counted with multiplicity, that are strictly less than 1.

In the case where $M^n$ is the critical catenoid, it is readily shown (c.f. Section 2.4) that $\mathrm{Ker}(\mathrm{J}^n)$ is one-dimensional, that $V$ is the space of functions with vanishing mean, and that the Steklov index of $M^n$ is related to its Morse index by the formula

$$\mathrm{MI}(n) = \mathrm{SI}(n) + 1. \tag{1.10}$$

We thus reformulate Theorem I in terms of the Steklov index as follows.

**Theorem II**

*For all $n$, let $SI(n)$ denote the Steklov index of the $n$-dimensional critical catenoid. Then*

$$\underset{n\rightarrow+\infty}{Lim} \frac{Log(SI(n))}{\sqrt{n}Log(\sqrt{n})} = 1. \tag{1.11}$$

**1.3 - Spherical harmonics.** We will use spherical harmonics to diagonalize the Dirichlet to Neumann operator (c.f. Section 2.3). To this end, let $\Delta^{S^{n-1}}$ denote the Laplace-Beltrami operator of the $(n-1)$-dimensional sphere $S^{n-1}$. The eigenvalues of this operator are

$$\lambda(n,m) := -m(n+m-2), \tag{1.12}$$

where $m$ here ranges over all non-negative integers. For each $m$, the space $\mathcal{H}(n,m)$ of $m$'th order *spherical harmonics* of dimension $(n-1)$ is defined to be the eigenspace of $\Delta^{S^{n-1}}$ with eigenvalue $\lambda(n,m)$. Since this space coincides with the space of restrictions to $S^{n-1}$ of $m$'th order homogeneous, harmonic polynomials over $\mathbb{R}^n$. In particular, it is straightforward to show that

$$\mathrm{Dim}(\mathcal{H}(n,m)) = \begin{pmatrix} n+m-1 \\ n-1 \end{pmatrix} - \begin{pmatrix} n+m-3 \\ n-1 \end{pmatrix}. \tag{1.13}$$

Now let $M^n$ be the critical catenoid discussed above. Up to rescaling, $\partial M^n$ is isometric to the disjoint union of two copies of $S^{n-1}$, so that functions over $\partial M^n$ identify with pair $(f,g)$ of functions over $S^{n-1}$. By symmetry, the eigenspaces of the Dirichlet to Neumann operator are then given by

$$\begin{aligned} \mathrm{W}_0(n,m) &:= \{(f,f)\ |\ f\in\mathcal{H}(n,m)\}\,,\ \text{and} \\ \mathrm{W}_1(n,m) &:= \{(f,-f)\ |\ f\in\mathcal{H}(n,m)\}\,, \end{aligned} \tag{1.14}$$

where $m$ again ranges over all non-negative integers. Note that DN is not actually defined over $\mathrm{W}_0(n,0)$, but we nonetheless include this space for notational simplicity. For all $(n,m)$ and for each $i$, let $\Lambda_i(n,m)$ denote the unique eigenvalue of DN over the space $\mathrm{W}_i(n,m)$ where, by convention, we set $\Lambda_0(n,0) := -\infty$. Bearing in mind (1.13), in order to prove Theorems I and II, it is sufficient to determine for which values of $m$ and $i$ the Steklov eigenvalue $\Lambda_i(n,m)$ is less than 1. A sufficiently accurate estimate is then derived from the following result, which we believe to be of independent interest.

**Theorem III**

*There exists $A>0$ such that for $1\leq m<n/A\mathrm{Log}(n)$,*

$$\begin{aligned} \Lambda_0(n,m) &= -\frac{n}{m}+m+O\left(\mathrm{Log}(n)+\frac{m^2\mathrm{Log}(n)}{n}\right),\ \textit{and} \\ \Lambda_1(n,m) &= m+O\left(\frac{m^2\mathrm{Log}(n)}{n}\right). \end{aligned} \tag{1.15}$$

This result follows immediately from Theorem 3.6.1, below, with $\alpha=2(n-1)$ and $\beta=m$.

**1.4 - Acknowledgements.** The second author was in part supported by the Simons Foundation (award #279968). The third author would like to thank Prof. Richard Schoen for helpful discussions.

## 2 - The geometric framework.

**2.1 - The basic geometry of free boundary minimal hypersurfaces.** Let $f := f(n)$ be a smooth, real-valued function defined over some interval $]a, b[$. Let $M := M(n)$ be the hypersurface of revolution in $\mathbb{R}^{n+1}$ generated by rotating the graph of $f$ about the $x$-axis. In the sequel, $f$ will be called the *profile* of $M$. We suppose that $M$ is minimal, and we first derive the differential condition that this imposes on $f$. For all $x$, the flux of $M$ through the hyperplane $\mathbb{R}^n \times \{x\}$ is given by

$$\mathcal{F}(x) := \frac{\Omega_{n-1} f(x)^{n-1}}{\sqrt{1 + f_x^2}}, \tag{2.1}$$

where $\Omega_{n-1}$ here denotes the volume of the unit $(n-1)$-dimensional sphere in $\mathbb{R}^n$, and the subscript $x$ here denotes differentiation in the $x$-direction. Since $M$ is minimal, this quantity is constant (c.f. [7]), and we thereby obtain the minimal surface equation for $f$, namely

$$f_x = \sqrt{Cf^{2(n-1)} - 1},$$

where $C$ is some positive real number. Upon rescaling, we may suppose that the neck radius $C^{-1/2(n-1)}$ of $M$ is equal to 1, so that

$$f_x = \sqrt{f^{2(n-1)} - 1}. \tag{2.2}$$

This is the equation that will be studied throughout the rest of this paper. In particular, we have the following simple formula for the second derivative of $f$.

$$f_{xx} = (n-1)f^{2n-3}. \tag{2.3}$$

We now suppose that $M$ is free boundary minimal, that is, it is contained in the closed unit ball and meets the unit sphere orthogonally along its boundary. In this case, its profile is the solution $f$ of (2.2) rescaled by a factor $R := R(n)$ which we now determine. Let $L := L(n)$ be such that $]-L, L[$ is the maximal interval over which $f$ is defined. It follows from (2.3) that $f$ is strictly convex over this interval, and it follows by maximality that $f_x$ tends to infinity as $x$ tends to $\pm L$. From this we deduce that there exists a unique point $W := W(n)$ such that $0 < W < L$ and

$$Wf_x(W) = f(W). \tag{2.4}$$

We readily verify that the scaling factor $R$ is given by

$$R^2 := W^2 + f(W)^2. \tag{2.5}$$

Finally, throughout the sequel, we will use the following explicit parametrisation for $M$.

$$F : [-W, W] \times S^{n-1} \to M; (x, \omega) \mapsto \left(\frac{x}{R}, \frac{f(x)\omega}{R}\right), \tag{2.6}$$

where $S^{n-1}$ here denotes the unit sphere in $\mathbb{R}^n$.

**2.2 - The stability spectrum and the Jacobi operator.** Recall from the introduction that the stability operator $Q : H^1(M) \oplus H^1(M) \to \mathbb{R}$ is defined to be the second variation of the volume of $M$ with respect to infinitesimal perturbations in the normal direction, and is given by

$$Q(\phi,\psi) := \int_M \langle \nabla\phi, \nabla\psi\rangle - \mathrm{Tr}\left(A^2\right)\phi\psi \mathrm{dVol}^M - \int_{\partial M} \phi\psi \mathrm{dVol}^{\partial M}. \tag{2.7}$$

Likewise, a function $\phi$ in $H^1(M)$ is said to be an *eigenfunction* of $Q$ with *eigenvalue* $\lambda$ whenever

$$Q(\phi,\psi) = \lambda \int_M \phi\psi \mathrm{dVol}^M, \tag{2.8}$$

for all $\psi$ in $H^1(M)$, and the *stability spectrum* of $M$, which we denote by $\mathrm{Spec}(M)$, is defined to be the set of all eigenvalues of $Q$ in $\mathbb{C}$. It follows from classical Sturm-Liouville theory (c.f. [8]) that the stability spectrum is discrete, every eigenfunction is smooth, every eigenvalue has finite multiplicity, and all but finitely many eigenvalues are strictly positive.

The stability spectrum may be characterised in other ways that are easier to study. Indeed, the *Jacobi operator* $\mathrm{J} := \mathrm{J}(n)$ and the *Robin boundary operator* $\mathrm{B} := \mathrm{B}(n)$ of $M$ are defined by

$$\begin{aligned} \mathrm{J}\phi &:= -\Delta^M\phi - \mathrm{Tr}\left(A^2\right)\phi, \text{ and} \\ \mathrm{B}\phi &:= \partial_\nu\phi - \phi\circ e, \end{aligned} \tag{2.9}$$

where $\Delta^M$ here denotes the Laplace-Beltrami operator of $M$, $\partial_\nu$ denotes the operator of differentiation in the outward-pointing unit, conormal, direction over $\partial M$ and $e : \partial M \to M$ denotes the canonical embedding. Using Stokes' theorem, we readily verify (c.f. [6], Section 5) that a smooth function $\phi$ is an eigenfunction of $Q$ with eigenvalue $\lambda$ whenever

$$\begin{aligned} \mathrm{J}\phi &= \lambda\phi, \text{ and} \\ \mathrm{B}\phi &= 0. \end{aligned}$$

In other words, the stability spectrum coincides with the spectrum of the Jacobi operator with Robin boundary conditions.

**Lemma 2.2.1**

*In the parametrisation (2.6), the Jacobi operator and the Robin boundary operator of M are given by*

$$\begin{aligned} J(n)\phi &= -\frac{R(n)^2}{f(n)^2}\Delta^{S^{n-1}}\phi - \frac{R(n)^2}{f(n)^{2(n-1)}}\phi_{xx} - \frac{n(n-1)R(n)^2}{f(n)^{2n}}\phi, \text{ and} \\ B_\pm(n)\phi &= \pm W(n)\phi_x(W(n)) - \phi(W(n)), \end{aligned} \tag{2.10}$$

*where $\Delta^{S^{n-1}}$ here denotes the Laplace-Beltrami operator of $S^{n-1}$, and $B_\pm(n)$ denotes the Robin boundary operator at $x = \pm W(n)$.*

**Proof:** Indeed, bearing in mind (2.2), the metric of $M$ is given in this parametrisation by

$$g = \frac{f^2}{R^2}g^{S^{n-1}} \oplus \frac{f^{2(n-1)}}{R^2}dx^2,$$

where $g^{S^{n-1}}$ here denotes the metric of $S^{n-1}$. In particular,

$$\mathrm{Det}(g) = \frac{f^{4(n-1)}}{R^{2(n+1)}},$$

and so the Laplacian is given by

$$\begin{aligned}\Delta\phi &= \frac{1}{\sqrt{\mathrm{Det}(g)}}\nabla_i\sqrt{\mathrm{Det}(g)}g^{ij}\nabla_j\phi\\ &= \frac{R^2}{f^2}\Delta^{S^{n-1}}\phi + \frac{R^2}{f^{2(n-1)}}\phi_{xx}.\end{aligned}$$

Next, since $M$ is a surface of revolution, its prinicipal curvature in every angular direction is

$$\kappa_\theta = \frac{R}{f^n},$$

and since it is minimal, its principal curvature in the longitudinal direction is

$$\kappa_l = -\frac{(n-1)R}{f^n}.$$

It follows that

$$\mathrm{Tr}\left(A^2\right) = (n-1)\kappa_\theta^2 + \kappa_l^2 = n(n-1)f^{-2n},$$

and the formula for the Jacobi operator follows. Finally, we readily verify that in this parametrisation, the vector $\pm\partial_x$ is an outward-pointing conormal to $\partial M$ in $M$ at $x = \pm W$ and that its length is equal to $1/W$. The formula for the Robin boundary operator follows, and this completes the proof. □

**2.3 - Jacobi fields and spherical harmonics.** A *Jacobi field* over $M$ is defined to be any function $\phi : M \to \mathbb{R}$ solving

$$\mathrm{J}\phi = 0.$$

We will study Jacobi fields using the classical technique of separation of variables which we now review. Consider first the Laplace-Beltrami operator $\Delta^{S^{n-1}}$ of $S^{n-1}$. The eigenvalues of this operator are given by (c.f. Chapter III.C.1 of [4])

$$\lambda(m,n) := -m(n+m-2),\ m \in \{0,1,2,...\}\,, \tag{2.11}$$

and, for each $m$, its eigenspace with eigenvalue $\lambda_m$ is the space $\mathcal{H}(n,m)$ of $(n-1)$-dimensional spherical harmonics of order $m$, which coincides with the space of restrictions to $S^{n-1}$ of homogeneous, harmonic polynomials of order $m$ over $\mathbb{R}^n$.

**Lemma 2.3.1**

*For all $n$ and for all $m$,*

$$\mathit{Dim}(\mathcal{H}(n,m)) = \begin{pmatrix} n+m-1 \\ n-1 \end{pmatrix} - \begin{pmatrix} n+m-3 \\ n-1 \end{pmatrix}. \tag{2.12}$$

**Proof:** Indeed, for all $(n,m)$, let $H(n,m)$ denote the space of homogeneous polynomials of order $m$ over $\mathbb{R}^n$. We readily verify that

$$\mathrm{Dim}(H(n,m)) = \begin{pmatrix} n+m-1 \\ n-1 \end{pmatrix}.$$

Now let $\Delta$ denote the standard Laplace-Beltrami operator of $\mathbb{R}^n$. Since this operator maps $H(n,m)$ surjectively into $H(n,m-2)$, the short sequence

$$0 \longrightarrow \mathcal{H}(n,m) \longrightarrow H(n,m) \xrightarrow{\Delta} H(n,m-2) \longrightarrow 0$$

is exact, and the result now follows. □

Since the Jacobi operator of $M$ is invariant under rotation, it preserves each harmonic mode and, for all $(n,m)$, we define $\mathrm{J}(n,m)$ to be the restriction of $\mathrm{J}(n)$ to the $m$'th harmonic mode. It follows immediately from (2.10) and (2.11) that

$$\mathrm{J}(n,m)\phi = -\frac{R(n)^2}{f(n)^{2(n-1)}}\phi_{xx} + \frac{R(n)^2 m(n+m-2)}{f(n)^2}\phi - \frac{R(n)^2 n(n-1)}{f^{2n}}\phi. \tag{2.13}$$

Consider now a Jacobi field $\phi$ over $M$. By separation of variables, $\phi$ may be written in the form

$$\phi(x,\omega) = \sum_{m=0}^{\infty} a_m(\omega)\phi_0(n,m)(x) + \sum_{m=0}^{\infty} b_m(\omega)\phi_1(n,m)(x), \tag{2.14}$$

where, for all $(n,m)$, the functions $a_m$ and $b_m$ are elements of $\mathcal{H}(n,m)$ and, for each $i$, $\phi_i(n,m)$ is the unique solution of

$$\mathrm{J}(n,m)\phi_i(n,m) = 0, \tag{2.15}$$

with initial conditions

$$\begin{aligned} \phi_i(n,m)(0) &= 1-i, \text{ and} \\ \phi_i(n,m)_x(0) &= i. \end{aligned} \tag{2.16}$$

In particular, for all $(n,m)$, $\phi_0(n,m)$ is even and $\phi_1(n,m)$ is odd. These functions will play a fundamental role in the sequel, and we will refer to them respectively as the *even and odd Jacobi fields* in the $m$'th *harmonic mode*.

Certain Jacobi fields can be determined explicitly, and the resulting formulae and the properties they imply will be of considerable use in what follows. Of most use to us will be the Jacobi fields in the zeroeth and first order harmonic modes, which arise from Killing fields of the ambient space.

**Lemma 2.3.2** *For all $n$,*

$$\begin{aligned}\phi_0(n,0) &= -xf(n)_x f(n)^{1-n} + f(n)^{2-n},\\ \phi_1(n,0) &= \frac{1}{(n-1)} f(n)_x f(n)^{1-n},\\ \phi_0(n,1) &= f(n)^{1-n}, \text{ and}\\ \phi_1(n,1) &= \frac{1}{n}(f(n)_x f(n)^{2-n} + xf(n)^{1-n}).\end{aligned} \tag{2.17}$$

**Proof:** Let $N : M \to S^n$ denote the outward-pointing unit normal vector field over $M$, so that, at the point $p := (x, f(x)\omega)$,

$$N(p) = \frac{1}{\sqrt{1+f_x^2}}(-f_x, \omega).$$

Let $X_{0,0}(x,y) := (x,y)$ denote the Killing field of infinitesimal dilatations of $\mathbb{R}^{n+1}$. Then

$$\phi_0(n,0) = \langle N, X_{0,0}\rangle = \frac{-xf_x + f}{\sqrt{1+f_x^2}} = -xf_x f^{1-n} + f^{2-n}.$$

Let $X_{0,1}(x,y) := (1,0)$ denote the Killing field of infinitesimal translations in the $x$-direction. Then

$$\phi_1(n,0) = \lambda\langle N, X_{0,1}\rangle = \lambda\frac{-f_x}{\sqrt{1+f_x^2}} = -\lambda f_x f^{1-n},$$

for some constant $\lambda$ which we readily verify is equal to $1/(n-1)$. Let $X_{1,0}(x,y) := (0,1)$ denote the Killing field of infinitesimal translations on the $y$-direction. Then

$$\phi_0(n,1) = \langle N, X_{1,0}\rangle = \frac{1}{\sqrt{1+f_x}} = f^{1-n}.$$

Finally, let $X_{1,1}(x,y) := (-y,x)$ denote the Killing field of infinitesimal rotations in the $(x-y)$-plane. Then

$$\phi_1(n,1) = \lambda\langle N, X_{1,1}\rangle = \lambda\frac{ff_x + x}{\sqrt{1+f_x^2}} = \lambda(f_x f^{2-n} + xf^{2-n}),$$

for some constant $\lambda$ which we readily verify is equal to $1/n$. This completes the proof. □

The following basic properties of Jacobi fields will prove useful.

**Lemma 2.3.3**

*For all $n$, $\phi_0(n,0)$ is strictly positive over $]-W(n),W(n)[$ and vanishes at $\pm W(n)$. In particular, the kernel of $J(n)$ over $M$ with Dirichlet boundary conditions is the linear subspace of $C^\infty(M)$ generated by this function.*

**Proof:** Indeed, denote $\phi := \phi_0(n,0)$. Bearing in mind (2.4), we have

$$\phi(W) = -Wf_x(W)f(W)^{1-n} + f(W)^{2-n} = 0.$$

Likewise, $\phi(-W) = 0$. Next, for all $x \in ]-W,W[$, by convexity,

$$f(x) > xf_x(x),$$

so that $\phi(x) > 0$. This proves the desired properties of $\phi$. The final assertion follows by the maximum principle, and this completes the proof. □

**Lemma 2.3.4**

*For all $(n,m)$ and for each $i$ such that $(m,i) \neq (0,0)$, $\phi_i(n,m)$ is strictly positive over $]0,W(n)]$.*

**Proof:** Consider first the case where $i = 0$. Denote

$$\sigma := \frac{\phi_0(n,m)_x}{\phi_0(n,m)}, \text{ and}$$
$$\tau := \frac{\phi_0(n,0)_x}{\phi_0(n,0)}.$$

By Lemma 2.3.3, $\tau$ is finite over $[0,W[$ but tends to $-\infty$ as $x$ tends to $\pm W$. At every point of $[0,W[$ where $\sigma$ is finite, we have

$$\sigma_x + \sigma^2 = m(n+m-2)f^{2(n-2)} - n(n-1)f^{-2} > -n(n-1)f^{-2} = \tau_x + \tau^2,$$

and since $\sigma(0) = 0 = \tau(0)$, it follows that $\sigma \geq \tau$, with equality only at 0. In particular, $\sigma$ is finite over $[0,W]$, and it follows that $\phi_0(n,m)$ never vanishes over this interval, as desired. The case where $i = 1$ is identical except that the initial condition on $\sigma$ becomes

$$\operatorname*{Lim}_{x\rightarrow 0^+} \sigma(x) = +\infty > \tau(0).$$

This completes the proof. □

Other special Jacobi fields are derived from the structure of the Jacobi operator. Although they will not be used much in the sequel, we include them for completeness and illustrative purposes.

**Lemma 2.3.5**

*For all $n$,*

$$\begin{aligned}\phi_0(n,n-1) &= \frac{(n-2)}{(n-1)}f(n)^{n-1}+\frac{1}{(n-1)}f(n)^{1-n},\\ \phi_0(n,n) &= f(n)^n,\\ \phi_1(n,2n-2) &= (3n-4)f(n)_xf(n)^{n-1}+f(n)_xf(n)^{1-n},\ \textit{and}\\ \phi_1(n,2n-1) &= f(n)_xf(n)^n.\end{aligned}\tag{2.18}$$

**Remark:** Intriguingly, by investigating the action of $\mathrm{J}(n,m)$ on the function $x^pf_x^qf^r$, for arbitrary $p$, $q$ and $r$, it may be possible to derive explicit formulae for all Jacobi fields in all harmonic modes. Although this will not be studied further in this paper, such an approach, were it to work, would yield an alternative proof of our main result.

**Proof:** Indeed, we readily verify that, for all $p$,

$$\begin{aligned}\mathrm{J}(n,m)f^p = &-\left(p(n+p-2)-m(n+m-2)\right)R^2f^{p-2}\\ &-\left(n(n-1)-p(p-1)\right)R^2f^{p-2n}.\end{aligned}$$

Likewise, for all $p$,

$$\begin{aligned}\mathrm{J}(n,m)f_xf^p = &-\left((p+n-1)(p+2n-3)-m(n+m-2)\right)R^2f_xf^{p-2}\\ &-\left(n(n-1)-p(p-1)\right)R^2f_xf^{p-2n}.\end{aligned}$$

The result follows upon combining these relations. $\square$

**2.4 - Steklov eigenvalues.** Using the parametrisation (2.6), square integrable functions over $\partial M$ are represented as pairs $(\psi_-,\psi_+)$ of functions in $L^2(S^{n-1})$. By classical Sturm-Liouville theory, every such function can be written in the form

$$(\psi_-,\psi_+)=\sum_{m=0}^{\infty}(a_m,a_m)+\sum_{m=0}^{\infty}(-b_m,b_m),\tag{2.19}$$

where, for all $m$, the functions $a_m$ and $b_m$ are elements of $\mathcal{H}(n,m)$. Now let $V$ be the orthogonal complement in $L^2(S^{n-1})\oplus L^2(S^{n-1})$ of the linear subspace generated by the pair $(1,1)$, that is

$$V:=\langle(1,1)\rangle^{\perp}.$$

By Lemma 2.3.3, $V$ is the space, constructed in Section 1.2 of the introduction, over which the Dirichet to Neumann operator is defined. Using (2.14), we readily show that if $(\psi_-,\psi_+)\in V$ is given by (2.19), then $\mathrm{DN}(\psi_-,\psi_+)$ is given explicitly by

$$\begin{aligned}\mathrm{DN}(\psi_-,\psi_+)=&\sum_{m=1}^{\infty}\frac{W(n)\phi_0(n,m)_x(W(n))}{\phi_0(n,m)(W(n))}(a_m,a_m)\\ &+\sum_{m=0}^{\infty}\frac{W(n)\phi_1(n,m)_x(W(n))}{\phi_1(n,m)(W(n))}(-b_m,b_m).\end{aligned}\tag{2.20}$$

In particular, the Steklov eigenvalues of $M$ are given by

$$\begin{aligned}\Lambda_0(n,m) &:= \frac{W(n)\phi_0(n,m)_x(W(n))}{\phi_0(n,m)(W(n))},\ m>0,\ \text{and}\\ \Lambda_1(n,m) &:= \frac{W(n)\phi_1(n,m)_x(W(n))}{\phi_1(n,m)(W(n))},\ m\geq 0,\end{aligned} \tag{2.21}$$

For completeness, we also define

$$\Lambda_0(n,0) := -\infty. \tag{2.22}$$

In the sequel, for all $(n,m)$, $\Lambda_0(n,m)$ and $\Lambda_1(n,m)$ will be called respectively the *even and odd Steklov eigenvalues* of $M$ in the $m$'th *harmonic mode*.

**Lemma 2.4.1**

*(1) for all $n$,*

$$\Lambda_0(n,n) = n \geq 1;$$

*(2) For all $n$,*

$$\Lambda_1(n,1) = 1;$$

*(3) for all $(n,m)$,*

$$\Lambda_1(n,m) > \Lambda_0(n,m);\ \textit{and}$$

*(4) for all $n$, for all $m > m'$ and for each $i$,*

$$\Lambda_i(n,m) > \Lambda_i(n,m').$$

**Proof:** To prove (1), recall that, by (2.18)

$$\phi_0(n,n) = f^n,$$

so that

$$\phi_0(n,n)_x = nf^{n-1}f_x,$$

and so, bearing in mind (2.4),

$$\Lambda_0(n,n) = \frac{W\phi_0(n,n)_x(W)}{\phi_0(n,n)(W)} = \frac{nf(W)^{n-1}Wf_x(W)}{f(W)^n} = n,$$

as desired.

To prove (2), recall that $\phi_1(n,1)$ is the Jacobi field of infinitesimal rotations. However, since rotations preserve the angle between $M$ and the unit sphere, we have

$$\mathrm{B}_\pm(n)\phi = 0,$$

and the result follows.

To prove (4) in the case where $i=0$, denote

$$\sigma := \frac{\phi_0(n,m)_x}{\phi_0(n,m)},\ \text{and}$$
$$\tau := \frac{\phi_0(n,m')_x}{\phi_0(n,m')}.$$

Then, as in the proof of Lemma 2.3.4, $\sigma \geq \tau$ with equality only at $x=0$. It follows that

$$\Lambda_0(n,m) = W\sigma(W) > W\tau(W) = \Lambda_0(n,m'),$$

as desired. Finally, (3) and (4) in the case where $i=1$ are verified in a similar manner, and this completes the proof. □

We conclude this section by showing how knowledge of the Steklov eigenvalues yields the Morse Index of $M$. First, for all $n$, and for each $i$, define $\mathrm{K}_i := \mathrm{K}_i(n)$ by

$$\mathrm{K}_i(n) = \#\left\{m \mid \Lambda_i(n,m) < 1\right\}. \tag{2.23}$$

By Lemma 2.4.1, (4), knowledge of $\mathrm{K}_0(n)$ and $\mathrm{K}_1(n)$ yields all harmonic modes with negative Steklov eigenvalues. Next, by (2.12), for all $(n,m)$ and for each $i$, the multiplicity of $\Lambda_i(n,m)$ is given by

$$\mathrm{Mult}(\Lambda_i(n,m)) = \mathrm{Dim}(\mathcal{H}(n,m)) = \begin{pmatrix} n+m-1 \\ n-1 \end{pmatrix} - \begin{pmatrix} n+m-3 \\ n-1 \end{pmatrix}, \tag{2.24}$$

so that this in turn yields the Steklov index of $M$. The Morse index of $M$ is then given by the following

**Lemma 2.4.2**

*For all $n$,*

$$\mathit{MI}(n) = \mathit{SI}(n) + 1. \tag{2.25}$$

**Proof:** First, by symmetry, the contributions to the Morse index of each of the harmonic modes can be studied independently of one another. Thus, fix $n$, $m$ and $i$. For $\lambda < 0$, let $\tilde{\phi}_i(n,m,\lambda)$ solve

$$\mathrm{J}(n,m)\tilde{\phi}_i(n,m,\lambda) = \lambda\tilde{\phi}_i(n,m,\lambda),$$

with initial conditions

$$\tilde{\phi}_i(n,m,\lambda)(0) = 1-i,\ \text{and}$$
$$\tilde{\phi}_i(n,m,\lambda)_x(0) = i,$$

and define

$$\tilde{\Lambda}_i(n,m,\lambda) = \frac{W(n)\tilde{\phi}_i(n,m,\lambda)_x(W(n))}{\tilde{\phi}_i(n,m,\lambda)(W(n))}.$$

Observe now that $\tilde{\phi}_i(n,m,\lambda)$ is an eigenfunction of $\mathrm{J}(n,m)$ for the Robin boundary problem if and only if $\tilde{\Lambda}_i(n,m,\lambda)=1$. However, reasoning as in the proof of Lemma 2.3.4, we see that $\tilde{\Lambda}_i(n,m,\lambda)$ is a strictly increasing function of $(-\lambda)$ which tends to infinity as $\lambda$ tends to minus infinity. It follows that there exists at most one value of $\lambda$ for which $\tilde{\phi}_i(n,m,\lambda)$ is an eigenfunction, and that this occurs if and only if

$$\Lambda_i(n,m)=\tilde{\Lambda}_i(n,m,0)<1.$$

The result follows. $\square$

Finally, we illustrate this technique by showing how the above information is already sufficient to recover the main result of [9], [24] and [21].

**Theorem 2.4.3, Devyver, Smith & Zhou, Tran**

*The Morse index of the critical catenoid in* $\mathbb{R}^3$ *is equal to* 4.

**Proof:** Indeed, by Lemma 2.4.1, for all $n$,

$$2\leq K_0(n)\leq n, \tag{2.26}$$

and

$$K_1(n)=1. \tag{2.27}$$

In particular, when $n=2$, $K_0(n)=2$. Recalling now that $\Lambda_0(n,0)=-\infty$ is not considered as a Steklov eigenvalue of $M$, it follows by (2.24) that

$$\mathrm{SI}(n)=\mathrm{Dim}(\mathcal{H}(2,0))+\mathrm{Dim}(\mathcal{H}(2,1))=3,$$

and it follows by (2.25) that

$$\mathrm{MI}(n)=1+\mathrm{SI}(n)=4,$$

as desired. $\square$

Observe that (2.26) and (2.27) are in fact valid for all $n$. Furthermore, when $n=3$, it follows from (2.18) that $\Lambda_0(n,2)>1$ if and only if

$$H(2):=f(2)(W(2))>3^{\frac{1}{4}}.$$

As this is readily verified by numerical means (c.f. Table 3, below), it follows that $\mathrm{MI}(3)=5$. However, as the dimension grows, the second inequality in (2.26) looses its utility, and the growth of the index is controlled by another, entirely different, phenomenon which we now proceed to study.

## 3 - The analytic solution.

**3.1 - The continuous family.** At this stage, it will be helpful to reformulate the problem to allow for non-integer values of the dimension. Thus, for all $\alpha \geq 2$, we define $f := f_\alpha$ to be the unique solution of

$$\begin{aligned} f_x &= \sqrt{f^\alpha - 1}, \text{ and} \\ f(0) &= 1. \end{aligned} \tag{3.1}$$

Let $L := L_\alpha$ be such that $]-L, L[$ is the maximal interval over which $f$ is defined, let $W := W_\alpha \in ]0, L[$ be the unique point such that

$$W f_x(W) = f(W), \tag{3.2}$$

and define $H := H_\alpha$ by

$$H := f(W). \tag{3.3}$$

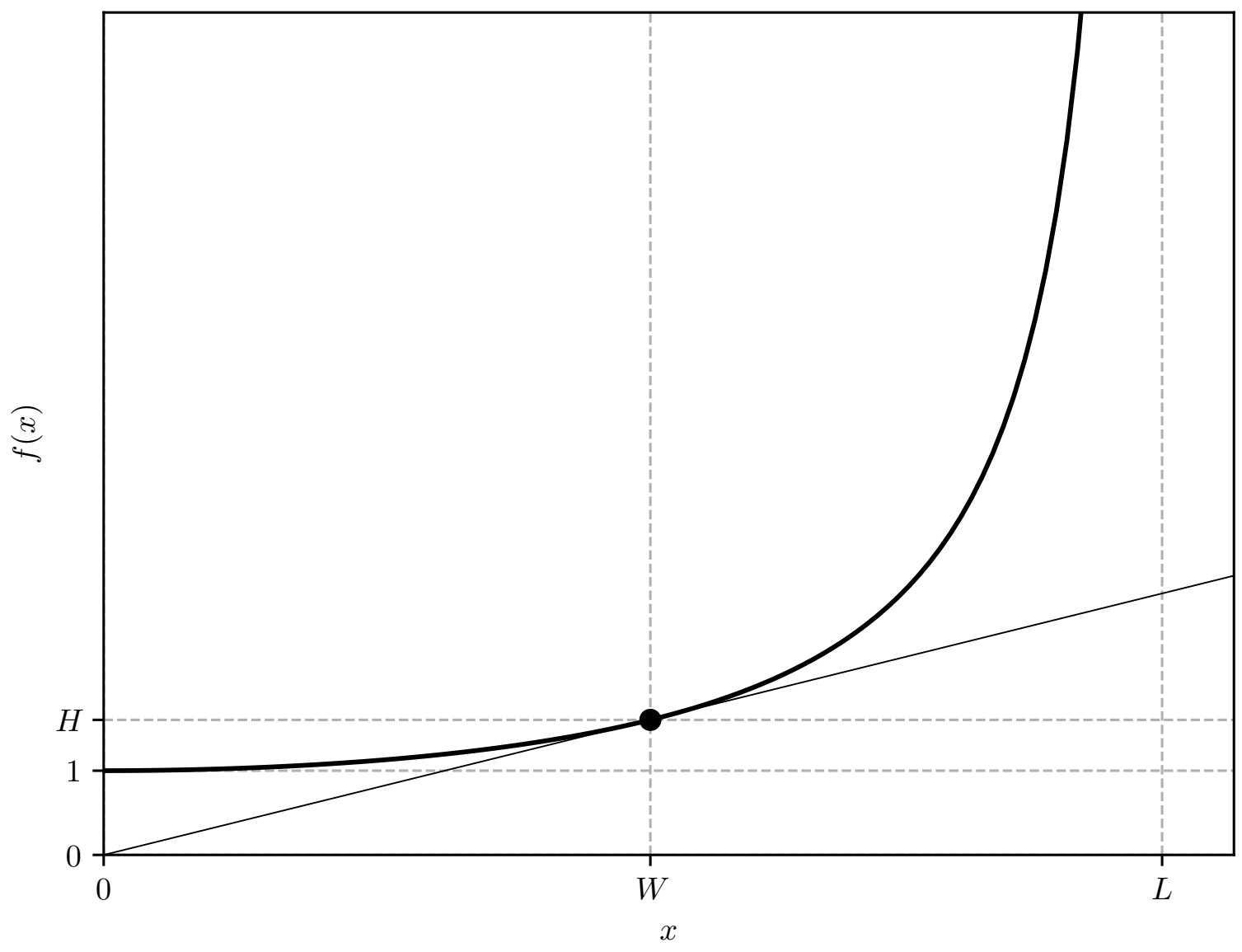

Figure 1 - The geometric meanings of the quantities $L$, $W$ and $H$ in relation to the graph of $f$.

The graph of $f$ together with the geometric meanings of the quantities $L$, $W$ and $H$ are illustrated in Figure 1. The remainder of this section is devoted to deriving asymptotic formulae for these quantities in terms of $\alpha$. We show

**Theorem 3.1.1**

*For all $\alpha$, let $L_\alpha$, $H_\alpha$ and $W_\alpha$ be defined as above. There exists a constant $\kappa > 0$ such that*

$$\begin{aligned} L_\alpha &= \frac{\pi}{\alpha}\left(1 + \frac{\kappa}{\alpha} + O\left(\frac{1}{\alpha^2}\right)\right), \\ H_\alpha &= \left(\frac{(\alpha-2)L_\alpha}{\alpha}\right)^{-\frac{2}{(\alpha-2)}}\left(1 + O\left(\frac{1}{\alpha^3}\right)\right), \textit{and} \\ W_\alpha &= L_\alpha\left(1 - \frac{2}{\alpha} + O\left(\frac{1}{\alpha^2}\right)\right). \end{aligned} \tag{3.4}$$

Theorem 3.1.1 is proven in Lemmas 3.1.2, 3.1.3 and 3.1.4, below. The constant $\kappa$, which will reappear throughout the sequel, is given by

$$\kappa := -\frac{1}{\pi}\beta_s\left(\frac{1}{2},\frac{1}{2}\right), \tag{3.5}$$

where $\beta$ here denotes Euler's Beta function and the subscript $s$ denotes partial differentiation with respect to the first variable. Recall that *Euler's Beta function* is given by

$$\beta(s,t) := \int_0^1 r^{s-1}(1-r)^{t-1}dr. \tag{3.6}$$

In particular, a straightforward integral yields

$$\beta\left(\frac{1}{2},\frac{1}{2}\right) = \pi. \tag{3.7}$$

In fact, upon expressing the derivative of $\beta$ in terms of the digamma function (c.f. [17]), we obtain

$$\beta_s\left(\frac{1}{2},\frac{1}{2}\right) = -2\pi\mathrm{Log}(2).$$

from which it follows that

$$\kappa = 2\mathrm{Log}(2).$$

However, we will have no need for this in the sequel.

Now let $g := g_\alpha : [1,\infty[\rightarrow [0,L]$ be the inverse function of $f$. Observe that

$$g(y) := \int_1^y \frac{1}{\sqrt{s^\alpha - 1}}ds. \tag{3.8}$$

In particular, for all $y \in ]1,\infty[$,

$$g_y(y) = \sum_{m=0}^{\infty}\binom{2m}{m}\frac{1}{4^m}y^{-\frac{(2m+1)\alpha}{2}}, \tag{3.9}$$

and integrating from $y$ to infinity yields, for $\alpha > 2$,

$$g(y) = L - \sum_{m=0}^{\infty}\binom{2m}{m}\frac{1}{4^m}\frac{2}{[(2m+1)\alpha-2]}y^{\frac{[2-(2m+1)\alpha]}{2}}. \tag{3.10}$$

**Lemma 3.1.2**

$$L_\alpha = \frac{\pi}{\alpha}\left(1 + \frac{\kappa}{\alpha} + O\left(\frac{1}{\alpha^2}\right)\right). \tag{3.11}$$

**Proof:** Indeed,

$$\begin{aligned} L &= \int_1^\infty \frac{1}{\sqrt{t^\alpha - 1}}dt \\ &= \frac{1}{\alpha}\beta\left(\frac{1}{2}-\frac{1}{\alpha},\frac{1}{2}\right), \end{aligned}$$

and the result follows by Taylor's Theorem together with (3.5) and (3.7). □

**Lemma 3.1.3**

$$H_\alpha = \left(\frac{(\alpha-2)L_\alpha}{\alpha}\right)^{-\frac{2}{(\alpha-2)}}\left(1+O\left(\frac{1}{\alpha^3}\right)\right). \tag{3.12}$$

**Proof:** Indeed, by (3.2) together with the inverse function theorem,

$$Hg_y(H) = g(H).$$

Substituting this into (3.9) and (3.10) yields

$$\sum_{m=0}^{\infty}\binom{2m}{m}\frac{1}{4^m}\frac{(2m+1)\alpha}{[(2m+1)\alpha-2]}H^{\frac{[2-(2m+1)\alpha]}{2}} = L,$$

so that

$$H^{-\alpha}\left(\sum_{m=0}^{\infty}\binom{2m}{m}\frac{1}{4^m}\frac{(2m+1)(\alpha-2)}{[(2m+1)\alpha-2]}H^{-m\alpha}\right)^{\frac{2\alpha}{(\alpha-2)}} = \left(\frac{(\alpha-2)L}{\alpha}\right)^{\frac{2\alpha}{(\alpha-2)}}.$$

Consider now the function

$$p(\omega,\zeta) := \sum_{m=0}^{\infty}\binom{2m}{m}\frac{1}{4^m}\frac{(2m+1)(1-2\omega)}{[(2m+1)-2\omega]}\zeta^m.$$

This function is analytic for $(\omega,\zeta)$ near $(0,0)$ and satisfies $p(\omega,0)=1$ for all $\omega$. It follows that

$$q(\omega,\zeta) := \zeta p(\omega,\zeta)^{\frac{2}{(1-2\omega)}}$$

is also analytic for $(\omega,\zeta)$ near $(0,0)$. Furthermore, $q(\omega,0)=0$ and $q_\zeta(\omega,0)=1$. It follows by the inverse function theorem for analytic functions that there exists a function $r(\omega,\xi)$ which is analytic in a neighbourhood of $(0,0)$ such that $q(\omega,r(\omega,\xi))=\xi$. Furthermore, $r(\omega,0)=0$ and $r_\xi(\omega,0)=1$. Since $L$ tends to zero as $\alpha$ tends to infinity, it follows that

$$H^{-\alpha} = r\left(\frac{1}{\alpha},\left(\frac{(\alpha-2)L}{\alpha}\right)^{\frac{2\alpha}{(\alpha-2)}}\right).$$

Finally, since

$$L = \mathrm{O}\left(\frac{1}{\alpha}\right),$$

it follows that

$$H^{-\alpha} = \left(\frac{(\alpha-2)L}{\alpha}\right)^{\frac{2\alpha}{(\alpha-2)}}\left(1+\mathrm{O}\left(\frac{1}{\alpha^2}\right)\right),$$

so that

$$H = \left(\frac{(\alpha-2)L}{\alpha}\right)^{-\frac{2}{(\alpha-2)}}\left(1+\mathrm{O}\left(\frac{1}{\alpha^3}\right)\right).$$

as desired. □

**Lemma 3.1.4**

$$W_\alpha = L_\alpha\left(1 - \frac{2}{\alpha} + O\left(\frac{1}{\alpha^2}\right)\right). \tag{3.13}$$

**Proof:** Indeed,

$$W = g(H) = Hg_y(H).$$

Substituting (3.9) and (3.12) into the above formula yields

$$W = \frac{(\alpha-2)L}{\alpha}\left(1 + \mathrm{O}\left(\frac{1}{\alpha^2}\right)\right),$$

and the result follows. □

**3.2 - The rescaled family.** For all $\alpha$, define the function $c := c_\alpha$ over the interval $[-1,1]$ by

$$c(x) := f(Lx)^{-\frac{(\alpha-2)}{2}}. \tag{3.14}$$

Observe that $c(0) = 1$, $c(\pm 1) = 0$, and

$$c_x(x) = -\frac{L(\alpha-2)}{2}\sqrt{1 - c(x)^{\frac{2\alpha}{(\alpha-2)}}}. \tag{3.15}$$

Define also the function $c_\infty$ over this interval by

$$c_\infty(x) := \cos\left(\frac{\pi x}{2}\right). \tag{3.16}$$

Observe that $c_\infty(0) = 1$, and

$$c_{\infty,x} = -\frac{\pi}{2}\sqrt{1 - c_\infty^2}. \tag{3.17}$$

Since $c$ is concave for all $\alpha \in [2,\infty]$, it is straightforward to show that $c_\alpha$ converges uniformly to $c_\infty$ as $\alpha$ tends to $+\infty$. However, we will require more refined estimates. In this section, we prove

**Theorem 3.2.1**

*For all $\alpha$, let $c_\alpha : [-1,1] \to \mathbb{R}$ be as in (3.14). There exists a constant $A > 0$ such that, for all $\alpha$,*

$$1 - \frac{A}{\alpha} \leq \frac{c_\alpha}{c_\infty} \leq 1 + \frac{A}{\alpha} \tag{3.18}$$

Theorem 3.2.1 is proven in Lemmas 3.2.2 and 3.2.5, below. First observe that, for large $\alpha$,

$$\frac{(\alpha-2)L}{\pi} = 1 - \frac{(2-\kappa)}{\alpha} + \mathrm{O}\left(\frac{1}{\alpha^2}\right) < 1. \tag{3.19}$$

**Lemma 3.2.2**

*For sufficiently large $\alpha$,*

$$\frac{c_\alpha}{c_\infty} \geq \frac{(\alpha-2)L_\alpha}{\pi} = 1 - \frac{(2-\kappa)}{\alpha} + O\left(\frac{1}{\alpha^2}\right). \tag{3.20}$$

**Proof:** Indeed, denote

$$\phi := \frac{(\alpha-2)L}{\pi}c_\infty.$$

Then, for sufficiently large $\alpha$,

$$\phi(0) = \frac{(\alpha-2)L}{\pi} = \left(1 - \frac{(2-\kappa)}{\alpha} + \mathrm{O}\left(\frac{1}{\alpha^2}\right)\right) < 1 = c(0).$$

On the other hand, $\phi(1) = 0 = c(1)$. Now suppose that $c - \phi$ attains a strictly negative minimum at some point $x \in ]0,1[$. Then, at this point,

$$0 < c(x) < \phi(x) < 1,$$

so that

$$c(x)^{\frac{2\alpha}{(\alpha-2)}} < c(x)^2 < \phi(x)^2 < c_\infty(x)^2.$$

In particular,

$$c_x(x) < \phi_x(x).$$

This is absurd, and the result follows. □

Upper bounds for $c$ require more work. First, we have

**Lemma 3.2.3**

*If $c_\alpha/c_\infty$ attains its maximal value at an interior point $x$ of the interval $]0,1[$, then*

$$\frac{c_\alpha}{c_\infty} \leq c_\infty(x)^{-\frac{2}{\alpha}}.$$

**Proof:** Indeed,

$$\begin{aligned} \mathrm{Log}(c)_x &= -\frac{(\alpha-2)L}{2}c^{-1}\sqrt{1-c^{\frac{2\alpha}{\alpha-2}}}, \text{ and} \\ \mathrm{Log}(c_\infty)_x &= -\frac{\pi}{2}c_\infty^{-1}\sqrt{1-c_\infty^2}. \end{aligned}$$

At the point $x$, these two quantities are equal so that, bearing in mind (3.20),

$$\frac{\sqrt{1-c_\infty^2}}{\sqrt{1-c^{\frac{2\alpha}{\alpha-2}}}} = \left(\frac{c_\infty}{c}\right)\frac{(\alpha-2)L}{\pi} \leq 1.$$

Consequently,

$$1 - c_\infty^2 \leq 1 - c^{\frac{2\alpha}{\alpha-2}},$$

and the result follows. □

**Lemma 3.2.4**

*There exists $\epsilon > 0$ and an analytic function $\Psi$ such that, for all $\alpha > 1/\epsilon$, and for all $1-\epsilon < x < 1$,*

$$c_\alpha^{\frac{\alpha}{(\alpha-2)}}(x) = \Psi\left(\frac{1}{\alpha}, \xi_\alpha(x)\right) sin\left(\xi_\alpha(x)\right),$$

*where*

$$\xi_\alpha(x) := \left(\frac{(\alpha-2)L_\alpha}{2}\right)^{\frac{\alpha}{(\alpha-2)}} (1-x)^{\frac{\alpha}{(\alpha-2)}}.$$

*Furthermore, for all $|\xi|, |\omega| < \epsilon$,*

$$\Psi(0,\xi) = \Psi(\omega, 0) = 1.$$

**Proof:** Observe first that if $y = f(Lx)$, then

$$x = \frac{g(y)}{L},$$

so that

$$(1-x) = \frac{1}{L}\sum_{m=0}^{\infty}\binom{2m}{m}\frac{1}{4^m}\frac{2}{[(2m+1)\alpha-2]}y^{\frac{[2-(2m+1)\alpha]}{2}}.$$

It follows that

$$\left(\frac{L(\alpha-2)}{2}(1-x)\right)^{\frac{\alpha}{(\alpha-2)}} = y^{-\frac{\alpha}{2}}\left(\sum_{m=0}^{\infty}\binom{2m}{m}\frac{1}{4^m}\frac{(\alpha-2)}{[(2m+1)\alpha-2]}y^{-m\alpha}\right)^{\frac{\alpha}{(\alpha-2)}}.$$

Consider now the function

$$p(\omega,\zeta) := \sum_{m=0}^{\infty}\binom{2m}{m}\frac{1}{4^m}\frac{(1-2\omega)}{(2m+1)-2\omega}\zeta^{2m}.$$

This function is analytic for $(\omega,\zeta)$ in a neighbourhood of $(0,0)$. Since $p(\omega,0)=1$, it follows that

$$q(\omega,\zeta) := \zeta p(\omega,\zeta)^{\frac{1}{(1-2\omega)}}$$

is also analytic for $(\omega,\zeta)$ in a neighbourhood of $(0,0)$, and

$$q_\zeta(\omega,0) = 1.$$

In particular, by the inverse function theorem for analytic functions, there exists an analytic function $r(\omega,\xi)$ such that, for all $(\omega,\xi)$ in a neighbourhood of $(0,0)$,

$$q(\omega, r(\omega,\xi)) = \xi.$$

Observe that $r(\omega,0)=0$ and $r_\xi(\omega,0)=1$, so that

$$r(\omega,\xi) = \Psi(\omega,\xi)\text{sin}(\xi),$$

for some function $\Psi$ which is also analytic near 0. However, since

$$q(0,\zeta) = \sum_{m=0}^{\infty}\binom{2m}{m}\frac{1}{4^m}\frac{1}{(2m+1)}\zeta^{2m+1} = \text{arcsin}(\zeta),$$

we have

$$r(0,\xi) = \sin(\xi),$$

and the result follows. □

**Lemma 3.2.5**

*There exists a constant $A>0$ such that, for all sufficiently large $\alpha$,*

$$\frac{c_\alpha}{c_\infty}\leq 1+\frac{A}{\alpha}. \tag{3.21}$$

**Proof:** Indeed, if $\epsilon$ is as in Lemma 3.2.4, then there exists $A>0$ such that, for all sufficiently large $\alpha$ and for all $1-\epsilon<x<1$,

$$\frac{c(x)}{c_\infty(x)}\leq 1+\frac{A}{\alpha}.$$

If $c/c_\infty$ attains its maximum over $[1-\epsilon,1]$, then we are done. Otherwise, $c/c_\infty$ attains its maximum over $[0,1-\epsilon]$. In this case, by Lemma 3.2.3,

$$\frac{c}{c_\infty}\leq \sin\left(\frac{\epsilon\pi}{2}\right)^{-\frac{2}{\alpha}}=1+\frac{2}{\alpha}\left|\mathrm{Log}\left(\frac{\epsilon\pi}{2}\right)\right|+\mathrm{O}\left(\frac{1}{\alpha^2}\right),$$

and the result follows. □

**3.3 - The rescaled Jacobi operators.** For all $\alpha\geq 2$ and for all $\beta\geq 0$, consider now the operator $\mathrm{J}(\alpha,\beta)$ given by

$$\mathrm{J}(\alpha,\beta)\phi:=\phi_{xx}-\left(A(\alpha,\beta)-B(\alpha)\right)\phi, \tag{3.22}$$

where

$$\begin{aligned}A(\alpha,\beta)&:=\frac{\beta(\alpha+2\beta-2)L_\alpha^2}{2c_\alpha^2},\ \text{and}\\ B(\alpha)&:=\frac{\alpha(\alpha+2)L_\alpha^2}{4}c_\alpha^{\frac{4}{(\alpha-2)}}.\end{aligned} \tag{3.23}$$

For each $i\in\{0,1\}$, define the function $\phi_i(\alpha,\beta):[-1,1]\to\mathbb{R}$ to be the unique solution of

$$\mathrm{J}(\alpha,\beta)\phi_i(\alpha,\beta)=0, \tag{3.24}$$

with initial conditions

$$\begin{aligned}\phi_i(\alpha,\beta)(0)&=1-i,\\ \phi_i(\alpha,\beta)_x(0)&=i.\end{aligned} \tag{3.25}$$

Observe that when $\alpha=2(n-1)$ and $\beta=m$, the operator $\mathrm{J}(\alpha,\beta)$ is - up to rescaling and multiplication by a smooth function - the $m$'th harmonic component of the Jacobi operator of the $n$-dimensional free boundary minimal catenoid in $\mathbb{R}^{n+1}$ which we studied in detail in Section 2.3. Via this observation, we derive explicit formulae for $\phi_i(\alpha,\beta)$ in certain special cases. Indeed, from Lemma 2.3.2, we readily obtain

**Theorem 3.3.1**

*For all $\alpha$, setting $t := L_\alpha x$,*

$$\begin{aligned}
\phi_0(\alpha,0) &= -t\sqrt{f_\alpha(t)^\alpha - 1}f_\alpha(t)^{-\frac{\alpha}{2}} + f_\alpha(t)^{\frac{(2-\alpha)}{2}},\\
\phi_1(\alpha,0) &= \frac{2}{\alpha L_\alpha}\sqrt{f_\alpha(t)^\alpha - 1}f_\alpha(t)^{-\frac{\alpha}{2}},\\
\phi_0(\alpha,1) &= f_\alpha^{-\frac{\alpha}{2}}(t),\ \text{and}\\
\phi_1(\alpha,1) &= \frac{2}{(2+\alpha)L_\alpha}\left(f_\alpha(t)^{\frac{(2-\alpha)}{2}}\sqrt{f_\alpha(t)^\alpha - 1} + tf_\alpha(t)^{-\frac{\alpha}{2}}\right).
\end{aligned} \tag{3.26}$$

In this section, we determine asymptotic formulae for the values of these functions and their derivatives at a certain critical point. Indeed, we prove

**Theorem 3.3.2**

*There exists a constant $A \geq 1$ such that for all $2 \leq \gamma \leq \alpha \leq Exp(\gamma/A)$,*

$$\begin{aligned}
\phi_0(\gamma,1)\left(\frac{W_\alpha}{L_\alpha}\right) &= \frac{\pi}{\alpha} + O\left(\frac{Log(\alpha)}{\gamma\alpha}\right),\\
\phi_1(\gamma,1)\left(\frac{W_\alpha}{L_\alpha}\right) &= \frac{2}{\pi} + O\left(\frac{Log(\alpha)}{\gamma}\right),\\
\phi_0(\gamma,1)_x\left(\frac{W_\alpha}{L_\alpha}\right) &= -\frac{\pi}{2} + O\left(\frac{Log(\alpha)}{\gamma}\right),\ \text{and}\\
\phi_1(\gamma,1)_x\left(\frac{W_\alpha}{L_\alpha}\right) &= \frac{2\alpha}{\pi\gamma} + O\left(\frac{\alpha Log(\alpha)}{\gamma^2}\right).
\end{aligned} \tag{3.27}$$

**Remark:** A straightforward extension of this result also yields uniform bounds - independent of $\gamma$ and $\alpha$ - of the functions $\phi_0(\gamma,1)$ and $\phi_1(\gamma,1)$ as well as their first derivatives over the interval $[0, W_\alpha/L_\alpha]$.

Theorem 3.3.2 is proven in Lemmas 3.3.6, 3.3.7, 3.3.8 and 3.3.9, below. We first derive some preliminary estimates.

**Lemma 3.3.3**

*For all $2 \leq \gamma \leq \alpha$,*

$$c_\gamma\left(\frac{W_\alpha}{L_\alpha}\right) = \frac{\pi}{\alpha} + O\left(\frac{1}{\gamma\alpha}\right). \tag{3.28}$$

**Proof:** Indeed, by (3.4) and (3.18),

$$\begin{aligned}
c_\gamma\left(\frac{W_\alpha}{L_\alpha}\right) &= \left(1 + \mathrm{O}\left(\frac{1}{\gamma}\right)\right)\cos\left(\frac{\pi W_\alpha}{2L_\alpha}\right)\\
&= \left(1 + \mathrm{O}\left(\frac{1}{\gamma}\right)\right)\sin\left(\frac{\pi}{\alpha} + \mathrm{O}\left(\frac{1}{\alpha^2}\right)\right)\\
&= \frac{\pi}{\alpha}\left(1 + \mathrm{O}\left(\frac{1}{\gamma}\right)\right)\left(1 + \mathrm{O}\left(\frac{1}{\alpha}\right)\right),
\end{aligned}$$

as desired. □

**Lemma 3.3.4**

*For all $\alpha > \gamma$ such that $Log(\alpha) \ll \gamma$,*

$$f_\gamma\left(\frac{L_\gamma W_\alpha}{L_\alpha}\right) = 1 + O\left(\frac{Log(\alpha)}{\gamma}\right). \tag{3.29}$$

**Proof:** Indeed, using (3.28) together with the definition of $c_\gamma$, we obtain

$$\begin{aligned}
f_\gamma\left(\frac{L_\gamma W_\alpha}{L_\alpha}\right) &= c_\gamma\left(\frac{W_\alpha}{L_\alpha}\right)^{-\frac{2}{(\gamma-2)}} \\
&= \left(\frac{\pi}{\alpha}\right)^{-\frac{2}{(\gamma-2)}}\left(1 + \mathrm{O}\left(\frac{1}{\gamma\alpha^2}\right)\right)^{-\frac{2}{(\gamma-2)}} \\
&= \left(1 + \mathrm{O}\left(\frac{\mathrm{Log}(\alpha)}{\gamma}\right)\right),
\end{aligned}$$

as desired. □

**Lemma 3.3.5**

*For all $\alpha > \gamma$ such that $Log(\alpha) \ll \gamma$,*

$$f_{\gamma,x}\left(\frac{L_\gamma W_\alpha}{L_\alpha}\right) = \frac{\alpha}{\pi} + O\left(\frac{\alpha Log(\alpha)}{\gamma}\right). \tag{3.30}$$

**Proof:** Indeed, using (3.18), (3.28) and (3.29) together with the definition of $c_\gamma$, we obtain

$$\begin{aligned}
f_{\gamma,x}\left(\frac{L_\gamma W_\alpha}{L_\alpha}\right) &= \left(f_\gamma\left(\frac{L_\gamma W_\alpha}{L_\alpha}\right)^\gamma - 1\right)^{\frac{1}{2}} \\
&= \left(c_\gamma\left(\frac{W_\alpha}{L_\alpha}\right)^{-2} f_\gamma\left(\frac{L_\gamma W_\alpha}{L_\alpha}\right)^2 - 1\right)^{\frac{1}{2}} \\
&= \left(\left(\frac{\pi}{\alpha} + \mathrm{O}\left(\frac{1}{\gamma\alpha}\right)\right)^{-2}\left(1 + \mathrm{O}\left(\frac{\mathrm{Log}(\alpha)}{\gamma}\right)\right) - 1\right)^{\frac{1}{2}} \\
&= \frac{\alpha}{\pi}\left(1 + \mathrm{O}\left(\frac{\mathrm{Log}(\alpha)}{\gamma}\right)\right),
\end{aligned}$$

as desired. □

**Lemma 3.3.6**

*For all $\alpha > \gamma$ such that $Log(\alpha) \ll \gamma$,*

$$\phi_0(\gamma,1)\left(\frac{W_\alpha}{L_\alpha}\right) = \frac{\pi}{\alpha} + O\left(\frac{Log(\alpha)}{\gamma\alpha}\right). \tag{3.31}$$

**Proof:** Indeed, by (3.26), (3.28) and (3.29),

$$\begin{aligned}\phi_0(\gamma,1)\left(\frac{W_\alpha}{L_\alpha}\right) &= f_\gamma\left(\frac{L_\gamma W_\alpha}{L_\alpha}\right)^{-1} c_\gamma\left(\frac{W_\alpha}{L_\alpha}\right)\\ &= \left(1 + \mathrm{O}\left(\frac{\mathrm{Log}(\alpha)}{\gamma}\right)\right)\left(\frac{\pi}{\alpha} + \mathrm{O}\left(\frac{1}{\gamma\alpha}\right)\right)\\ &= \frac{\pi}{\alpha} + \mathrm{O}\left(\frac{\mathrm{Log}(\alpha)}{\gamma\alpha}\right),\end{aligned}$$

as desired. □

**Lemma 3.3.7**

*For all $\alpha > \gamma$ such that $Log(\alpha) \ll \gamma$,*

$$\phi_1(\gamma,1)\left(\frac{W_\alpha}{L_\alpha}\right) = \frac{2}{\pi} + O\left(\frac{Log(\alpha)}{\gamma}\right). \tag{3.32}$$

**Proof:** Indeed, by (3.4), (3.28) and (3.30),

$$\frac{2}{(2+\gamma)L_\gamma} c_\gamma\left(\frac{W_\alpha}{L_\alpha}\right) f_{\gamma,x}\left(\frac{L_\gamma W_\alpha}{L_\alpha}\right) = \frac{2}{\pi}\left(1 + \mathrm{O}\left(\frac{\mathrm{Log}(\alpha)}{\gamma}\right)\right).$$

Next, since $c_\gamma \leq 1$ and $f_\gamma \geq 1$, by (3.4) again,

$$\frac{2}{(2+\gamma)L_\gamma}\left(\frac{L_\gamma W_\alpha}{L_\alpha}\right) f_\gamma\left(\frac{L_\gamma W_\alpha}{L_\alpha}\right)^{-\frac{\alpha}{2}} = \mathrm{O}\left(\frac{1}{\gamma}\right).$$

The result now follows by (3.26). □

**Lemma 3.3.8**

*For all $\alpha > \gamma$ such that $Log(\alpha) \ll \gamma$,*

$$\phi_0(\gamma,1)_x\left(\frac{W_\alpha}{L_\alpha}\right) = -\frac{\pi}{2} + O\left(\frac{Log(\alpha)}{\gamma}\right). \tag{3.33}$$

**Proof:** Indeed, differentiating (3.26) yields, for all $x$,

$$\phi_0(\gamma,1)_x = -\frac{\gamma L_\gamma}{2} c_\gamma(x) f_\gamma^{-2}(L_\gamma x) f_{\gamma,x}(L_\gamma x).$$

Thus, by (3.4), (3.28), (3.29) and (3.30),

$$\begin{aligned}\phi_0(\gamma,1)_x\left(\frac{W_\alpha}{L_\alpha}\right) &= -\frac{\pi}{2}\left(1+\mathrm{O}\left(\frac{1}{\gamma}\right)\right)\left(\frac{\pi}{\alpha}+\mathrm{O}\left(\frac{1}{\gamma\alpha}\right)\right)\\ &\qquad\times\left(1+\mathrm{O}\left(\frac{\mathrm{Log}(\alpha)}{\gamma}\right)\right)\left(\frac{\alpha}{\pi}+\mathrm{O}\left(\frac{\alpha\mathrm{Log}(\alpha)}{\gamma}\right)\right)\\ &= -\frac{\pi}{2}+\mathrm{O}\left(\frac{\mathrm{Log}(\alpha)}{\gamma}\right),\end{aligned}$$

as desired. □

**Lemma 3.3.9**

*For all $\alpha > \gamma$ such that $Log(\alpha) \ll \gamma$,*

$$\phi_1(\gamma,1)_x\left(\frac{W_\alpha}{L_\alpha}\right) = \frac{2\alpha}{\pi\gamma}+O\left(\frac{\alpha Log(\alpha)}{\gamma^2}\right). \tag{3.34}$$

**Proof:** Indeed, differentiating (3.26) yields, for all $\gamma$,

$$\begin{aligned}\phi_1(\gamma,1)(x) = \frac{2}{(2+\gamma)}\bigg(&c_\gamma^{-1}(x)f_\gamma(L_\gamma x)+\frac{\gamma}{2}c_\gamma(x)f_\gamma(L_\gamma x)^{-1}\\ &-\frac{\gamma}{2}(L_\gamma x)c_\gamma(x)f_\gamma(L_\gamma x)^{-2}f_{\gamma,x}(L_\gamma x)\bigg).\end{aligned}$$

Thus, by (3.4), (3.28), (3.29) and (3.30),

$$\begin{aligned}\phi_1(\gamma,1)_x\left(\frac{W_\alpha}{L_\alpha}\right) &= \frac{2}{\gamma}\left(1+\mathrm{O}\left(\frac{1}{\gamma}\right)\right)\bigg[\frac{\alpha}{\pi}\left(1+\mathrm{O}\left(\frac{1}{\gamma}\right)\right)\left(1+\mathrm{O}\left(\frac{\mathrm{Log}(\alpha)}{\gamma}\right)\right)\\ &\qquad+\frac{\pi\gamma}{2\alpha}\left(1+\mathrm{O}\left(\frac{1}{\gamma}\right)\right)\left(1+\mathrm{O}\left(\frac{\mathrm{Log}(\alpha)}{\gamma}\right)\right)\\ &\qquad-\frac{\pi}{2}\left(1+\mathrm{O}\left(\frac{1}{\alpha}\right)\right)\left(1+\mathrm{O}\left(\frac{1}{\gamma}\right)\right)\\ &\qquad\qquad\times\left(1+\mathrm{O}\left(\frac{\mathrm{Log}(\alpha)}{\gamma}\right)\right)\left(1+\mathrm{O}\left(\frac{\mathrm{Log}(\alpha)}{\gamma}\right)\right)\bigg].\\ &= \frac{2\alpha}{\pi\gamma}\left(1+\mathrm{O}\left(\frac{1}{\gamma}\right)\right)\left(1+\mathrm{O}\left(\frac{\mathrm{Log}(\alpha)}{\gamma}\right)\right)\\ &= \frac{2\alpha}{\pi\gamma}+\mathrm{O}\left(\frac{\alpha\mathrm{Log}(\alpha)}{\gamma^2}\right),\end{aligned}$$

as desired. □

**3.4 - Elementary perturbation theory.** In what follows, we will use solutions of perturbations of the operator $\mathrm{J}(\gamma,1)$. To this end, in the present section, we recall some general theory. Consider therefore a second order, linear operator of the form

$$\mathcal{L}\phi := \phi'' + g\phi,$$

defined over some interval $[0,A]$, where $g:[0,A]\rightarrow\mathbb{R}$ is some continuous function. Let $a$ and $b$ be the unique solutions of the equation $\mathcal{L}\phi=0$ with initial conditions $a(0)=b'(0)=1$ and $a'(0)=b(0)=0$. Consider the kernel $K$ defined by

$$K(x,y) := b(x)a(y) - a(x)b(y), \tag{3.35}$$

and denote also by $K:L^1([0,A])\rightarrow C^0([0,A])$ the operator given by

$$Ku(x) := \int_0^x K(x,y)u(y)dy. \tag{3.36}$$

Observe that since $C^0([0,A])$ embeds continuously into $L^1([0,A])$, the operator $K$ can also be thought of as defining a bounded linear map from this space to itself.

**Lemma 3.4.1**

*For all $u\in C^0([0,A])$, the function $Ku$ is twice continuously differentiable. Furthermore,*

$$\mathcal{L}Ku = u,$$

*and*

$$Ku(0) = (Ku)_x(0) = 0.$$

**Proof:** Trivially $Ku$ is continuous and $Ku(0)=0$. Since $K(x,x)=0$, differentiating under the integral yields

$$(Ku)_x(x) = \int_0^x K_x(x,y)u(y)dy,$$

so that $Ku$ is continuously differentiable and $(Ku)_x(0)=0$. Now, differentiating a second time, and denoting the Wronskian of $a$ and $b$ by $W$, we obtain

$$\begin{aligned}(Ku)_{xx}(x) &= K_x(x,x)u(x) + \int_0^x K_{xx}(x,y)u(y)dy\\ &= W(x)u(x) + \int_0^x \left(b_{xx}(x)a(y) - a_{xx}(x)b(y)\right)u(y)dy\\ &= W(x)u(x) - g(x)\int_0^x K(x,y)u(y)dy\\ &= W(x)u(x) - g(x)(Ku)(x).\end{aligned}$$

It follows that $Ku$ is twice continuously differentiable and, since the Wronskian is constant and equal to 1,

$$\mathcal{L}Ku = u,$$

as desired. □

Consider a smooth function $h : [0, A] \to \mathbb{R}$ and denote by $M_h : C^0([0, A]) \to L^1([0, A])$ the operator of multiplication by this function. We now introduce the operator

$$\hat{K} := \sum_{m=0}^{\infty} (-KM_h)^m . \tag{3.37}$$

Since $K$ is an operator of Volterra type, this series converges for all $h$. Indeed, we have

**Theorem 3.4.2**

*The series* (3.37) *converges absolutely. In particular,* $\hat{K}$ *defines a bounded linear map from* $C^0([0, 1])$ *to itself and its operator norm satisfies*

$$\|\hat{K}\| \leq e^{\|K\|_{L^\infty}\|h\|_{L^1}} . \tag{3.38}$$

**Proof:** For all $m$, and for all $\phi \in C^0([0, A])$,

$$\begin{aligned}(KM_h)^m\phi(x) = \int_0^x K(x, x_1)h(x_1) \int_0^{x_1} K(x_1, x_2)h(x_2)...\\ \int_0^{x_{m-1}} K(x_{m-1}, x_m)h(x_m)\phi(x_m)dx_m...dx_1.\end{aligned}$$

Hence

$$\begin{aligned}|(KM_h)^m\phi(x)| &\leq \|K\|_{L^\infty}^m \|\phi\|_{L^\infty} \int_{0<x_m<...<x_1<x} |h(x_1)| ... |h(x_m)| \, dx_1...dx_m\\ &= \frac{1}{m!}\|K\|_{L^\infty}^m \|\phi\|_{L^\infty} \int_{0<x_1,...,x_m<x} |h(x_1)| ... |h(x_m)| \, dx_1...dx_m\\ &\leq \frac{1}{m!}\|K\|_{L^\infty}^m \|h\|_{L^1}^m \|\phi\|_{L^\infty} .\end{aligned}$$

The operator norm of $(KM_h)^m$ therefore satisfies

$$\|(KM_h)^m\| \leq \frac{1}{m!}\|K\|_{L^\infty}^m \|h\|_{L^1}^m,$$

and the result follows upon summing over all $m$. □

Consider now the functions $\hat{a}$ and $\hat{b}$ given by

$$\begin{aligned}\hat{a} &:= \hat{K}a, \text{ and}\\ \hat{b} &:= \hat{K}b.\end{aligned} \tag{3.39}$$

The main result of this section is

**Theorem 3.4.3**

*The functions $\hat{a}$ and $\hat{b}$ are the unique solutions of*

$$(\mathcal{L}+h)\phi = 0$$

*with initial conditions $\hat{a}(0)=\hat{b}_x(0)=1$ and $\hat{a}_x(0)=\hat{b}(0)=0$. Furthermore*

$$\begin{aligned}\hat{a}-a&=-KM_h\hat{a},\\ \hat{a}_x-a_x&=-K_xM_h\hat{a},\\ \hat{b}-b&=-KM_h\hat{b},\ \textit{and}\\ \hat{b}_x-b_x&=-K_xM_h\hat{b},\end{aligned}\tag{3.40}$$

*where $K_x$ here denotes the operator with kernel $K_x(x,y)$.*

**Proof:** It suffices to prove the result for $\hat{a}$, as the proof for $\hat{b}$ is identical. However,

$$\begin{aligned}\hat{a}&=a-KM_h\hat{a},\ \text{and}\\ \hat{a}_x&=a_x-K_xM_h\hat{a},\end{aligned}$$

and (3.40) follows. Furthermore, by Lemma 3.4.1, $\hat{a}$ satisfies the desired initial conditions. Finally, applying $\mathcal{L}$ to each side of the first relation in (3.40) yields

$$\mathcal{L}\hat{a}=\mathcal{L}a-\mathcal{L}KM_h\hat{a}=-h\hat{a},$$

and this completes the proof. □

**3.5 - Perturbations of the Jacobi operator.** Continuing to use the notation of the previous section, for $2\leq\gamma\leq\alpha$ such that $\mathrm{Log}(\alpha)\ll\gamma$, we set

$$\begin{aligned}\mathcal{L}&:=\mathrm{J}(\gamma,1),\ \text{and}\\ h&:=B(\alpha)-B(\gamma),\end{aligned}\tag{3.41}$$

and, for each $i\in\{0,1\}$, we define

$$\psi_i(\alpha,\gamma):=\hat{K}\phi_i(\gamma,1).\tag{3.42}$$

In this section, we determine the values of the functions $\psi_i(\alpha,\gamma)$ and their derivatives at the point $W_\alpha/L_\alpha$. We prove

**Theorem 3.5.1**

*With $A$ as in Theorem 3.3.2, for all $2 \leq \gamma \leq \alpha \leq Exp(\gamma/A)$,*

$$\begin{aligned}
\psi_0(\alpha,\gamma)\left(\frac{W_\alpha}{L_\alpha}\right) &= \frac{\pi}{\gamma} + O\left(\frac{Log(\alpha)}{\gamma^2}\right),\\
\psi_1(\alpha,\gamma)\left(\frac{W_\alpha}{L_\alpha}\right) &= \frac{2}{\pi} + O\left(\frac{Log(\alpha)}{\gamma}\right),\\
\psi_0(\alpha,\gamma)_x\left(\frac{W_\alpha}{L_\alpha}\right) &= -\frac{\pi}{2} + \frac{\pi(\alpha-\gamma)}{\gamma^2} + O\left(\frac{Log(\alpha)}{\gamma} + \frac{\alpha Log(\alpha)}{\gamma^3}\right), \text{ and}\\
\psi_1(\alpha,\gamma)_x\left(\frac{W_\alpha}{L_\alpha}\right) &= \frac{2\alpha}{\pi\gamma} + O\left(\frac{\alpha Log(\alpha)}{\gamma^2}\right).
\end{aligned} \tag{3.43}$$

**Proof:** This result follows from Theorem 3.3.2 and Lemmas 3.5.3, 3.5.5 and 3.5.6, below. We shall only present details of the proof for $\psi_0(\alpha,\gamma)_x$, as the remaining cases are similar. First, by (3.45), over the interval $[0, W_\alpha/L_\alpha]$,

$$\psi_0(\alpha,\gamma) = \phi_0(\alpha,\gamma) + \epsilon,$$

where

$$\epsilon = \mathrm{O}\left(\frac{1}{\gamma}\right).$$

Substituting this into (3.40) yields, over this interval

$$\psi_0(\alpha,\gamma)_x = \phi_0(\gamma,1)_x - K_x M_h \psi_0(\gamma,1) - K_x M_h \epsilon.$$

However, by (3.44),

$$\|M_h\| = \|h\|_{L^1} = \mathrm{O}\left(\frac{1}{\gamma}\right).$$

Likewise, by (3.27) and the definition of $K$,

$$K_x = \mathrm{O}\left(\frac{\alpha}{\gamma}\right).$$

It follows that

$$K_x M_h \epsilon = \mathrm{O}\left(\frac{\alpha}{\gamma^3}\right),$$

and the result now follows by (3.27), (3.47) and (3.48). □

We first estimate the difference between $\psi_i(\alpha,\gamma)$ and $\phi_i(\gamma,1)$.

**Lemma 3.5.2**

*For all* $2 \leq \gamma \leq \alpha$*,*

$$\|B(\alpha) - B(\gamma)\|_{L^1([0,W_\gamma/L_\gamma])} = O\left(\frac{1}{\gamma}\right). \tag{3.44}$$

**Proof:** First, using (3.4), we have

$$\begin{aligned}\int_0^{W_\alpha/L_\alpha}\left(1-\cos\left(\frac{\pi x}{2}\right)^{\frac{4}{(\gamma-2)}}\right)dx &= \int_0^{1-\frac{\alpha}{2}}\left(1-\cos\left(\frac{\pi x}{2}\right)^{\frac{4}{(\gamma-2)}}\right)dx + \mathrm{O}\left(\frac{1}{\alpha^2}\right)\\ &\leq \int_{\frac{2}{\alpha}}^1 1 - x^{\frac{4}{(\gamma-2)}}dx + \mathrm{O}\left(\frac{1}{\alpha^2}\right)\\ &= \frac{4}{\gamma} + \mathrm{O}\left(\frac{1}{\alpha\gamma}\right).\end{aligned}$$

Consequently, by (3.18) and (3.4), again,

$$\begin{aligned}B(\alpha) - B(\gamma) &= \frac{\pi^2}{4}\cos\left(\frac{\pi x}{2}\right)^{\frac{4}{(\gamma-2)}} - \frac{\pi^2}{4}\cos\left(\frac{\pi x}{2}\right)^{\frac{4}{(\alpha-2)}} + \mathrm{O}\left(\frac{1}{\gamma}\right)\\ &= \frac{\pi^2}{4}\left(1-\cos\left(\frac{\pi x}{2}\right)^{\frac{4}{(\alpha-2)}}\right) - \frac{\pi^2}{4}\left(1-\cos\left(\frac{\pi x}{2}\right)^{\frac{4}{(\gamma-2)}}\right) + \mathrm{O}\left(\frac{1}{\gamma}\right),\end{aligned}$$

and the result follows. □

**Lemma 3.5.3**

*For all* $2 \leq \gamma \leq \alpha$*, and for each* $i$*, over the interval* $[0, W_\alpha/L_\alpha]$*,*

$$\psi_i(\alpha,\gamma) = \phi_i(\gamma,1) + O\left(\frac{1}{\gamma}\right). \tag{3.45}$$

**Proof:** Indeed, by (3.40),

$$\psi_i(\alpha,\gamma) - \phi_i(\gamma,1) = -KM_h\psi_i(\alpha,\gamma).$$

However, by (3.44),

$$\|M_h\| = \|h\|_{L^1} = \mathrm{O}\left(\frac{1}{\gamma}\right).$$

Likewise, by (3.27) and the definition of $K$,

$$\|K\| = \mathrm{O}\left(1\right),$$

and the result follows. □

It remains to estimate the different components of $KM_h\phi_i$ and $K_xM_h\phi_i$.

**Lemma 3.5.4**

*For all* $2 \leq \gamma \leq \alpha$,

$$\int_0^{W_\alpha/L_\alpha} Log(c_\gamma)c_\gamma^{\frac{2(\gamma+2)}{(\gamma-2)}}\,dx = \frac{(1-\kappa)}{4} + O\left(\frac{1}{\gamma}\right), \tag{3.46}$$

*where* $\kappa$ *is the constant defined in* (3.5)*.*

**Proof:** Indeed, by (3.4) and (3.18),

$$\int_0^{W_\alpha/L_\alpha} \mathrm{Log}(c_\gamma)c_\gamma^{\frac{2(\gamma+2)}{(\gamma-2)}}\,dx = \int_0^1 \mathrm{Log}\left(\cos\left(\frac{\pi x}{2}\right)\right)\cos^2\left(\frac{\pi x}{2}\right)dx + \mathrm{O}\left(\frac{1}{\gamma}\right).$$

However,

$$\begin{aligned}\int_0^1 \mathrm{Log}\left(\cos\left(\frac{\pi x}{2}\right)\right)\cos^2\left(\frac{\pi x}{2}\right)dx &= \frac{1}{2\pi}\int_0^1 \mathrm{Log}(r)r^{\frac{1}{2}}(1-r)^{-\frac{1}{2}}dr\\ &= \frac{1}{2\pi}\int_0^1 \mathrm{Log}(r)r^{-\frac{1}{2}}(1-r)^{\frac{1}{2}}dr + \frac{1}{\pi}\int_0^1 r^{-\frac{1}{2}}(1-r)^{\frac{1}{2}}dr\\ &= \frac{1}{2\pi}\int_0^1 \mathrm{Log}(r)r^{-\frac{1}{2}}(1-r)^{-\frac{1}{2}}dr - \frac{1}{2\pi}\int_0^1 \mathrm{Log}(r)r^{\frac{1}{2}}(1-r)^{-\frac{1}{2}}dr\\ &\qquad + \frac{1}{\pi}\int_0^1 r^{-\frac{1}{2}}(1-r)^{\frac{1}{2}}dr,\end{aligned}$$

where the second line is obtained by integrating by parts. Consequently,

$$\begin{aligned}\int_0^1 \mathrm{Log}\left(\cos\left(\frac{\pi x}{2}\right)\right)\cos^2\left(\frac{\pi x}{2}\right)dx &= \frac{1}{4\pi}\int_0^1 \mathrm{Log}(r)r^{-\frac{1}{2}}(1-r)^{-\frac{1}{2}}dr\\ &\qquad + \frac{1}{2\pi}\int_0^1 r^{-\frac{1}{2}}(1-r)^{\frac{1}{2}}dr\\ &= \frac{1}{4\pi}\left(\beta_s\left(\frac{1}{2},\frac{1}{2}\right) + \beta\left(\frac{1}{2},\frac{1}{2}\right)\right),\end{aligned}$$

where $\beta$ here denotes Euler's beta function, and the subscript $s$ denotes partial differentiation with respect to the first component. The result now follows by definition of $\kappa$. $\square$

**Lemma 3.5.5**

*For all* $2 \leq \gamma \leq \alpha$,

$$\int_0^{W_\alpha/L_\alpha} \phi_0(\gamma,1)^2(B(\alpha)-B(\gamma))dx = \frac{\pi^2}{2\alpha} - \frac{\pi^2}{2\gamma} + O\left(\frac{1}{\gamma^2}\right). \tag{3.47}$$

**Proof:** First observe that,

$$\begin{aligned}\phi_0(\gamma,1)^2(B(\alpha)-B(\gamma)) &= f_\gamma(L_\gamma x)^{-\gamma}\left(\frac{\alpha(\alpha+2)L_\alpha^2}{4}c_\alpha^{\frac{4}{(\alpha-2)}}-\frac{\gamma(\gamma+2)L_\gamma^2}{4}c_\gamma^{\frac{4}{(\gamma-2)}}\right)\\ &= \left(\frac{\alpha(\alpha+2)L_\alpha^2}{4}-\frac{\gamma(\gamma+2)L_\gamma^2}{4}\right)c_\gamma^{\frac{2(\gamma+2)}{(\gamma-2)}}\\ &\qquad+\frac{\alpha(\alpha+2)L_\alpha^2}{4}c_\gamma^{\frac{2\gamma}{(\gamma-2)}}\left(c_\alpha^{\frac{4}{(\alpha-2)}}-c_\gamma^{\frac{4}{(\gamma-2)}}\right).\end{aligned}$$

However, by (3.4),

$$\left(\frac{\alpha(\alpha+2)L_\alpha^2}{4}-\frac{\gamma(\gamma+2)L_\gamma^2}{4}\right)=\frac{2(\kappa+1)\pi^2}{4}\left(\frac{1}{\alpha}-\frac{1}{\gamma}\right)+\mathrm{O}\left(\frac{1}{\gamma^2}\right),$$

so that, by (3.18),

$$\begin{aligned}&\int_0^{W_\alpha/L_\alpha}\left(\frac{\alpha(\alpha+2)L_\alpha^2}{4}-\frac{\gamma(\gamma+2)L_\gamma^2}{4}\right)c_\gamma^{\frac{2(\gamma+2)}{(\gamma-2)}}dx\\ &\qquad=\int_0^1\frac{2(\kappa+1)\pi^2}{4}\left(\frac{1}{\alpha}-\frac{1}{\gamma}\right)\cos^2\left(\frac{\pi x}{2}\right)dx+\mathrm{O}\left(\frac{1}{\gamma^2}\right)\\ &\qquad=\frac{(\kappa+1)\pi^2}{4}\left(\frac{1}{\alpha}-\frac{1}{\gamma}\right)+\mathrm{O}\left(\frac{1}{\gamma^2}\right).\end{aligned}$$

By (3.18) again,

$$\begin{aligned}c_\alpha^{\frac{4}{(\alpha-2)}}&=c_\gamma^{\frac{4}{(\alpha-2)}}\left(1+\mathrm{O}\left(\frac{1}{\gamma}\right)\right)^{\frac{4}{(\alpha-2)}}\\ &=c_\gamma^{\frac{4}{(\alpha-2)}}+\mathrm{O}\left(\frac{1}{\gamma\alpha}\right),\end{aligned}$$

so that

$$c_\gamma^{\frac{2\gamma}{(\gamma-2)}}\left(c_\alpha^{\frac{4}{(\alpha-2)}}-c_\gamma^{\frac{4}{(\gamma-2)}}\right)=c_\gamma^{\frac{2(\gamma+2)}{(\gamma-2)}+\left(\frac{4}{(\alpha-2)}-\frac{4}{(\gamma-2)}\right)}-c_\gamma^{\frac{2(\gamma+2)}{(\gamma-2)}}+\mathrm{O}\left(\frac{1}{\gamma\alpha}\right).$$

Differentiating under the integral and using (3.46), we obtain

$$\begin{aligned}&\int_0^{W_\alpha/L_\alpha}c_\gamma^{\frac{2\gamma}{(\gamma-2)}}\left(c_\alpha^{\frac{4}{(\alpha-2)}}-c_\gamma^{\frac{4}{(\gamma-2)}}\right)dx\\ &\qquad=\left(\frac{4}{(\alpha-2)}-\frac{4}{(\gamma-2)}\right)\int_0^{W_\alpha/L_\alpha}\mathrm{Log}(c_\gamma)c_\gamma^{\frac{2(\gamma+2)}{(\gamma-2)}}dx+\mathrm{O}\left(\frac{1}{\gamma^2}\right)\\ &\qquad=\left(\frac{1-\kappa}{\alpha}-\frac{1-\kappa}{\gamma}\right)+\mathrm{O}\left(\frac{1}{\gamma^2}\right).\end{aligned}$$

The result now follows upon combining these relations. □

**Lemma 3.5.6**

*For all* $2 \leq \gamma \leq \alpha$,

$$\begin{aligned}\int_0^{W_\alpha/L_\alpha} \phi_0(\gamma,1)\phi_1(\gamma,1)(B(\alpha)-B(\gamma))dx = O\left(\frac{1}{\gamma}\right), \text{ and}\\ \int_0^{W_\alpha/L_\alpha} \phi_1(\gamma,1)^2(B(\alpha)-B(\gamma))dx = O\left(\frac{1}{\gamma}\right).\end{aligned} \tag{3.48}$$

**Proof:** Indeed, it follows from (3.27) that $\phi_0(\gamma,1)$ and $\phi_1(\gamma,1)$ are uniformly bounded independent of $\gamma$, and the result now follows by (3.44). □

**3.6 - Steklov eigenvalues.** For all $\alpha \geq 2$, for all $\beta \geq 0$, and for each $i \in \{0,1\}$, define

$$\Lambda_i(\alpha,\beta) := \frac{W_\alpha \phi_i(\alpha,\beta)_x(W_\alpha/L_\alpha)}{L_\alpha\phi_i(\alpha,\beta)(W_\alpha/L_\alpha)}. \tag{3.49}$$

Observe that, when $\alpha = 2(n-1)$ and $\beta = m$, $\Lambda_0(\alpha,\beta)$ and $\Lambda_1(\alpha,\beta)$ are respectively the even and odd Steklov eigenvalues in the $m$'th harmonic mode of the $n$-dimensional free boundary minimal catenoid in $\mathbb{R}^{n+1}$, which we studied in detail in Section 2.4. In this section, we derive asymptotic formulae for these functions. Indeed, we prove

**Theorem 3.6.1**

*There exists* $A > 0$ *such that for* $1 \leq \beta < \alpha/ALog(\alpha)$,

$$\begin{aligned}\Lambda_0(\alpha,\beta) &= -\frac{\alpha}{2\beta} + \beta + O\left(Log(\alpha) + \frac{\beta^2 Log(\alpha)}{\alpha}\right), \text{ and}\\ \Lambda_1(\alpha,\beta) &= \beta + O\left(\frac{\beta^2 Log(\alpha)}{\alpha}\right).\end{aligned} \tag{3.50}$$

**Proof:** This follows immediately from Lemma 3.1.4 and Lemma 3.6.3 below. □

We begin by estimating the first term in the zeroeth order component of the operator $\mathrm{J}(\alpha,\beta)$.

**Lemma 3.6.2**

*There exists* $C > 0$ *such that if* $\alpha > C\beta$, *then*

$$A(\gamma_-,1) \leq A(\alpha,\beta) \leq A(\gamma_+,1),$$

*where*

$$\gamma_\pm := \frac{\alpha}{\beta} \mp 2C.$$

**Proof:** It suffices to prove the second inequality, as the proof of the first is almost identical. However, by (3.4), setting $\gamma := \gamma_+$, for sufficiently large $C$, we have

$$\begin{aligned}\frac{\beta(\alpha+2\beta-2)L_\alpha^2}{2c_\alpha^2} &= \frac{\pi^2}{2}\left(\frac{\beta}{\alpha}\right)\left(1+2\left(\frac{\beta}{\alpha}\right)-\frac{2}{\alpha}\right)\left(1+\mathrm{O}\left(\frac{1}{\alpha}\right)\right)c_\alpha^{-2}\\ &\leq \frac{\pi^2}{2\gamma}\left(1-\frac{C}{\gamma}+\mathrm{O}\left(\frac{1}{\gamma}\right)\right)\cos\left(\frac{\pi x}{2}\right)^{-2}.\end{aligned}$$

Likewise,

$$\frac{\gamma L_\gamma^2}{2c_\gamma^2} \geq \frac{\pi^2}{2\gamma}\left(1+\mathrm{O}\left(\frac{1}{\gamma}\right)\right)\cos\left(\frac{\pi x}{2}\right)^{-2},$$

and the result follows. $\square$

**Lemma 3.6.3**

*There exists $A>0$ such that if $\alpha > A$ and if $\beta \leq \alpha/ALog(\alpha)$, then*

$$\begin{aligned}\frac{\psi_0(\alpha,\beta)_x(W_\alpha/L_\alpha)}{\psi_0(\alpha,\beta)(W_\alpha/L_\alpha)} &= -\frac{\alpha}{2\beta}+\beta+O\left(Log(\alpha)+\frac{\beta^2 Log(\alpha)}{\alpha}\right),\ \text{and}\\ \frac{\psi_1(\alpha,\beta)_x(W_\alpha/L_\alpha)}{\psi_1(\alpha,\beta)(W_\alpha/L_\alpha)} &= \beta+O\left(\frac{\beta^2 Log(\alpha)}{\alpha}\right).\end{aligned}\tag{3.51}$$

**Proof:** We prove an asymptotic upper bound for the first quotient, as a lower bound for this quotient as well bounds for the second quotient are proven in the same manner. Thus, with $\gamma_\pm$ as in Lemma 3.6.2, define

$$\begin{aligned}\sigma &:= \frac{\phi_0(\alpha,\beta)_x}{\phi_0(\alpha,\beta)},\ \text{and}\\ \tau &:= \frac{\psi_0(\alpha,\gamma_+)_x}{\psi_0(\alpha,\gamma_+)}.\end{aligned}$$

As in Lemma 2.3.4, the functions $\phi_0(\alpha,\beta)$ and $\psi_0(\alpha,\gamma_+)$ are non-vanishing over the interval $[0, W_\alpha/L_\alpha]$ so that both $\sigma$ and $\tau$ are finite over this interval. Thus, bearing in mind the definition of $\gamma_+$, over this interval

$$\sigma_x+\sigma^2 = A(\alpha,\beta)-B(\alpha) \leq A(\gamma_+,1)-B(\alpha) = \tau_x+\tau^2.$$

Since $\sigma(0)=0=\tau(0)$, it follows that $\sigma\leq\tau$ at every point of $[0,W_\alpha/L_\alpha]$, and the result now follows by (3.43). $\square$

We now prove the main results of this paper. First recall the function $K_0(n)$ introduced in Section 2.4.

**Theorem 3.6.4**

*For all $n$,*

$$K_0(n) = \sqrt{n} + O\left(Log(n)\right). \tag{3.52}$$

**Proof:** First, for all $\alpha \geq 2$ define

$$\mathrm{K}(\alpha) := \mathrm{Sup}\left\{\beta \mid \Lambda_0(\alpha,\beta) < 1\right\},$$

and observe that

$$|K_0(n) - K(\alpha/2+1)| \leq 1,$$

so that it is sufficient to show that

$$K(\alpha) = \sqrt{\frac{\alpha}{2}} + \mathrm{O}\left(\mathrm{Log}(\alpha)\right).$$

Consider now the function

$$f(x) := x - \frac{1}{x}.$$

Observe that $f(1) = 0$ and that $f$ has non-vanishing derivative at this point. In particular, there exists $\delta, \epsilon > 0$ such that the restriction of $f$ to $[1-\delta, 1+\delta]$ is a diffeomorphism onto its image and $[-\epsilon, \epsilon] \subseteq f([-\delta, \delta])$. Thus, bearing in mind (3.50), for sufficiently large $\alpha$, there exists $x \in [1-\delta, 1+\delta]$ such that if $\beta := \sqrt{\alpha/2}x$, then

$$\Lambda_0(\alpha,\beta) = 1.$$

Furthermore, as in Lemma 2.3.4, this is the unique value of $\beta$ for which this holds, and it follows that

$$K_0(\alpha) = \beta.$$

Finally, observe that

$$f(x) = x - \frac{1}{x} = \mathrm{O}\left(\frac{\mathrm{Log}(\alpha)}{\sqrt{\alpha}}\right),$$

so that, by the inverse function theorem,

$$x = 1 + \mathrm{O}\left(\frac{\mathrm{Log}(\alpha)}{\sqrt{\alpha}}\right),$$

and this completes the proof. $\square$

**Theorem 3.6.5**

*For all $n$, let MI($n$) denote the Morse index of the $n$-dimensional free boundary minimal catenoid in $\mathbb{R}^{n+1}$. MI($n$) satisfies the following asymptotic relation as $n$ tends to infinity.*

$$\mathit{Log}(\mathit{MI}(n)) = \sqrt{n}\mathit{Log}(\sqrt{n}) + \sqrt{n} + O\left(\mathit{Log}(\sqrt{n})^2\right). \tag{3.53}$$

**Proof:** Indeed, bearing in mind (2.24), (2.27) and (2.25),

$$\mathrm{MI}(n) = 1 + \begin{pmatrix} n+k-2 \\ k-1 \end{pmatrix} + \begin{pmatrix} n+k-3 \\ k-2 \end{pmatrix},$$

where $k = \mathrm{K}_0(n)$. However, for $m \ll n$, we have

$$\begin{aligned} \mathrm{Log}\left(\begin{pmatrix} n+m \\ m-1 \end{pmatrix}\right) &= \int_{n+1}^{n+m} \mathrm{Log}(x)dx - \int_1^{m-1} \mathrm{Log}(x)dx + \mathrm{O}\left(\mathrm{Log}(n)\right) \\ &= (n+m)\mathrm{Log}(n+m) - (n+1)\mathrm{Log}(n+1) \\ &\qquad - (m-1)\mathrm{Log}(m-1) + \mathrm{O}\left(\mathrm{Log}(n)\right) \\ &= (n+m)\mathrm{Log}(n) + m - n\mathrm{Log}(n) - m\mathrm{Log}(m) + \mathrm{O}\left(\mathrm{Log}(n)\right) \\ &= m\mathrm{Log}\left(\frac{n}{m}\right) + m + \mathrm{O}\left(\mathrm{Log}(n)\right), \end{aligned}$$

and the result now follows by (3.52). $\square$

## 4 - Numerical results.

**4.1 - Overview.** Much of our work was inspired by the output of numerical experiments. Indeed, had the numerical data not revealed such intriguing patterns, we would not have expected such comprehensive results. For this reason, we conclude this paper with a detailed description of the numerical calculations that we used. Above all, we discuss methods for computing numerical values of the Steklov eigenvalue $\Lambda := \Lambda_0(n,m)$ for arbitrary $n$ and $m$. First recall from (2.21) that

$$\Lambda = \frac{W\phi_x(W)}{\phi(W)}, \tag{4.1}$$

where $W := W(n)$ and $\phi := \phi_0(n,m)$ is the solution of

$$\mathrm{J}(n,m)\phi = 0, \quad \phi(0) = 1, \quad \phi_x(0) = 0. \tag{4.2}$$

These values will be computed numerically via the following general, two-step procedure.

| $n$ | $W$ | $n$ | $W$ | $n$ | $W$ | $n$ | $W$ | $n$ | $W$ |
|---|---|---|---|---|---|---|---|---|---|
| | | 21 | 0.07733 | 41 | 0.03897 | 61 | 0.02605 | 81 | 0.01956 |
| 2 | 1.19968 | 22 | 0.07370 | 42 | 0.03802 | 62 | 0.02562 | 82 | 0.01932 |
| 3 | 0.67715 | 23 | 0.07040 | 43 | 0.03713 | 63 | 0.02521 | 83 | 0.01908 |
| 4 | 0.47300 | 24 | 0.06738 | 44 | 0.03627 | 64 | 0.02481 | 84 | 0.01886 |
| 5 | 0.36357 | 25 | 0.06461 | 45 | 0.03545 | 65 | 0.02443 | 85 | 0.01863 |
| 6 | 0.29529 | 26 | 0.06206 | 46 | 0.03467 | 66 | 0.02405 | 86 | 0.01841 |
| 7 | 0.24860 | 27 | 0.05970 | 47 | 0.03392 | 67 | 0.02369 | 87 | 0.01820 |
| 8 | 0.21465 | 28 | 0.05752 | 48 | 0.03320 | 68 | 0.02334 | 88 | 0.01799 |
| 9 | 0.18887 | 29 | 0.05548 | 49 | 0.03252 | 69 | 0.02300 | 89 | 0.01779 |
| 10 | 0.16861 | 30 | 0.05359 | 50 | 0.03186 | 70 | 0.02266 | 90 | 0.01759 |
| 11 | 0.15227 | 31 | 0.05182 | 51 | 0.03122 | 71 | 0.02234 | 91 | 0.01739 |
| 12 | 0.13882 | 32 | 0.05017 | 52 | 0.03061 | 72 | 0.02203 | 92 | 0.01720 |
| 13 | 0.12755 | 33 | 0.04862 | 53 | 0.03003 | 73 | 0.02172 | 93 | 0.01702 |
| 14 | 0.11798 | 34 | 0.04716 | 54 | 0.02947 | 74 | 0.02143 | 94 | 0.01683 |
| 15 | 0.10974 | 35 | 0.04578 | 55 | 0.02892 | 75 | 0.02114 | 95 | 0.01666 |
| 16 | 0.10258 | 36 | 0.04449 | 56 | 0.02840 | 76 | 0.02086 | 96 | 0.01648 |
| 17 | 0.09629 | 37 | 0.04326 | 57 | 0.02790 | 77 | 0.02058 | 97 | 0.01631 |
| 18 | 0.09073 | 38 | 0.04210 | 58 | 0.02741 | 78 | 0.02032 | 98 | 0.01614 |
| 19 | 0.08578 | 39 | 0.04100 | 59 | 0.02694 | 79 | 0.02006 | 99 | 0.01598 |
| 20 | 0.08134 | 40 | 0.03996 | 60 | 0.02648 | 80 | 0.01981 | 100 | 0.01582 |

Table 2 - Numerical values of $W$ to five decimal places for $n$=2,$\cdots$,100.

(1) Using numerical root-finding methods, we find $W > 0$ such that

$$\frac{Wf_x(W)}{f(W)} = 1. \tag{4.3}$$

In order to set up this root-finding problem, it will also be necessary to numerically calculate either $f(x)$ or its inverse $g(y)$ (c.f. Section 3.1). This is carried out respectively using numerical integration and numerical quadrature.

(2) Using numerical integration methods, we solve the initial value problem for $\phi$ over the interval $[0, W]$ and we compute $\Lambda = W\phi_x(W)/\phi(W)$.

In the following subsections, these two steps will be described in greater detail and the numerical values of various quantities will be presented for various values of $m$ and $n$.

**4.2 - Numerical estimates of $W$ and $H$.** We first find $W > 0$ solving (4.3). To do so, we interpret $W$ as the unique root of the function

$$x \mapsto \frac{xf_x(x)}{f(x)} - 1, \tag{4.4}$$

over the range $x > 0$, to which a standard numerical root-finding method can then be applied. In order to compute (4.4), it is necessary to determine numerically both $f(x)$ and

| $n$ | $H$ | $n$ | $H$ | $n$ | $H$ | $n$ | $H$ | $n$ | $H$ |
|---|---|---|---|---|---|---|---|---|---|
| | | 21 | 1.14435 | 41 | 1.08678 | 61 | 1.06379 | 81 | 1.05106 |
| 2 | 1.81017 | 22 | 1.13939 | 42 | 1.08519 | 62 | 1.06298 | 82 | 1.05057 |
| 3 | 1.60312 | 23 | 1.13479 | 43 | 1.08366 | 63 | 1.06220 | 83 | 1.05009 |
| 4 | 1.48937 | 24 | 1.13053 | 44 | 1.08219 | 64 | 1.06144 | 84 | 1.04962 |
| 5 | 1.41609 | 25 | 1.12657 | 45 | 1.08077 | 65 | 1.06070 | 85 | 1.04916 |
| 6 | 1.36434 | 26 | 1.12286 | 46 | 1.07942 | 66 | 1.05998 | 86 | 1.04871 |
| 7 | 1.32557 | 27 | 1.11940 | 47 | 1.07811 | 67 | 1.05927 | 87 | 1.04826 |
| 8 | 1.29526 | 28 | 1.11615 | 48 | 1.07684 | 68 | 1.05859 | 88 | 1.04783 |
| 9 | 1.27082 | 29 | 1.11310 | 49 | 1.07563 | 69 | 1.05792 | 89 | 1.04740 |
| 10 | 1.25063 | 30 | 1.11022 | 50 | 1.07445 | 70 | 1.05727 | 90 | 1.04699 |
| 11 | 1.23363 | 31 | 1.10750 | 51 | 1.07332 | 71 | 1.05664 | 91 | 1.04658 |
| 12 | 1.21909 | 32 | 1.10493 | 52 | 1.07222 | 72 | 1.05602 | 92 | 1.04618 |
| 13 | 1.20648 | 33 | 1.10249 | 53 | 1.07116 | 73 | 1.05542 | 93 | 1.04578 |
| 14 | 1.19543 | 34 | 1.10018 | 54 | 1.07014 | 74 | 1.05483 | 94 | 1.04540 |
| 15 | 1.18566 | 35 | 1.09799 | 55 | 1.06914 | 75 | 1.05425 | 95 | 1.04502 |
| 16 | 1.17695 | 36 | 1.09589 | 56 | 1.06818 | 76 | 1.05369 | 96 | 1.04465 |
| 17 | 1.16912 | 37 | 1.09390 | 57 | 1.06725 | 77 | 1.05314 | 97 | 1.04428 |
| 18 | 1.16205 | 38 | 1.09200 | 58 | 1.06634 | 78 | 1.05261 | 98 | 1.04392 |
| 19 | 1.15562 | 39 | 1.09018 | 59 | 1.06547 | 79 | 1.05208 | 99 | 1.04357 |
| 20 | 1.14974 | 40 | 1.08845 | 60 | 1.06461 | 80 | 1.05157 | 100 | 1.04322 |

Table 3 - Numerical values of $H$ to five decimal places for $n$=2,$\cdots$,100.

$f_x(x)$ for arbitrary $x > 0$. To do so, we apply a numerical integrator to the initial value problem

$$f_{xx} = \frac{(n-1)(1+f_x^2)}{f}, \quad f(0) = 1, \quad f_x(0) = 0.$$

over the interval $[0, W]$. We use the Dormand and Prince 8(5,3) Runge–Kutta method (specifically, the DOP853 routine described in [13]) with adaptive step-size control to compute $f(x)$ and $f_x(x)$ to within a desired error tolerance. We then apply Brent's root-finding method [5] to (4.4) and this yields $W$.

Another approach is to transform the root-finding problem for $W$ into an equivalent problem for $H = f(W)$. Indeed, bearing in mind (3.8), we readily verify that $H$ is the unique root of the function

$$y \mapsto \frac{g(y)\sqrt{y^{2(n-1)}-1}}{y} - 1. \tag{4.5}$$

To compute (4.5) for arbitrary $H > 1$, we require numerical values of the integral

$$g(y) = \int_1^y \frac{1}{\sqrt{s^{2(n-1)}-1}} ds$$

| $n \backslash m$ | 2 | 3 | 4 | 5 | 6 | 7 | 8 | 9 | 10 |
|---|---|---|---|---|---|---|---|---|---|
| 2 | 2.00000 | 3.42087 | 4.68211 | 5.91042 | 7.12720 | 8.33838 | 9.54634 | 10.75228 | 11.95684 |
| 3 | 1.47402 | 3.00000 | 4.21410 | 5.35118 | 6.46346 | 7.56556 | 8.66253 | 9.75649 | 10.84851 |
| 4 | 1.03016 | 2.73719 | 4.00000 | 5.13565 | 6.22693 | 7.30018 | 8.36494 | 9.42512 | 10.48256 |
| 5 | 0.59915 | 2.50035 | 3.83842 | 5.00000 | 6.09515 | 7.16209 | 8.21573 | 9.26246 | 10.30522 |
| 6 | 0.16852 | 2.26598 | 3.68957 | 4.88937 | 6.00000 | 7.07098 | 8.12279 | 9.16458 | 10.20075 |
| 7 | −0.26493 | 2.02810 | 3.54180 | 4.78632 | 5.91901 | 7.00000 | 8.05523 | 9.09677 | 10.13062 |
| 8 | −0.70195 | 1.78520 | 3.39115 | 4.68416 | 5.84302 | 6.93791 | 8.00000 | 9.04432 | 10.07851 |
| 9 | −1.14258 | 1.53710 | 3.23629 | 4.58010 | 5.76788 | 6.87935 | 7.95076 | 9.00000 | 10.03642 |
| 10 | −1.58662 | 1.28405 | 3.07690 | 4.47298 | 5.69157 | 6.82158 | 7.90413 | 8.95992 | 10.00000 |
| 11 | −2.03378 | 1.02647 | 2.91305 | 4.36235 | 5.61311 | 6.76312 | 7.85822 | 8.92185 | 9.96671 |
| 12 | −2.48378 | 0.76476 | 2.74498 | 4.24812 | 5.53204 | 6.70316 | 7.81192 | 8.88441 | 9.93498 |
| 13 | −2.93635 | 0.49930 | 2.57298 | 4.13037 | 5.44817 | 6.64128 | 7.76459 | 8.84678 | 9.90382 |
| 14 | −3.39125 | 0.23046 | 2.39734 | 4.00925 | 5.36148 | 6.57725 | 7.71583 | 8.80843 | 9.87259 |
| 15 | −3.84828 | −0.04147 | 2.21836 | 3.88495 | 5.27201 | 6.51100 | 7.66543 | 8.76902 | 9.84086 |
| 16 | −4.30723 | −0.31619 | 2.03630 | 3.75765 | 5.17987 | 6.44249 | 7.61327 | 8.72836 | 9.80835 |
| 17 | −4.76796 | −0.59348 | 1.85140 | 3.62754 | 5.08515 | 6.37178 | 7.55931 | 8.68631 | 9.77486 |
| 18 | −5.23030 | −0.87311 | 1.66387 | 3.49482 | 4.98801 | 6.29891 | 7.50354 | 8.64280 | 9.74027 |
| 19 | −5.69414 | −1.15490 | 1.47392 | 3.35965 | 4.88854 | 6.22396 | 7.44598 | 8.59782 | 9.70452 |
| 20 | −6.15935 | −1.43867 | 1.28172 | 3.22219 | 4.78689 | 6.14702 | 7.38666 | 8.55135 | 9.66755 |

Table 4 - Numerical values of $\Lambda$ to five decimal places for $m$=2,$\cdots$,10 and $n$= 2,$\cdots$,20.

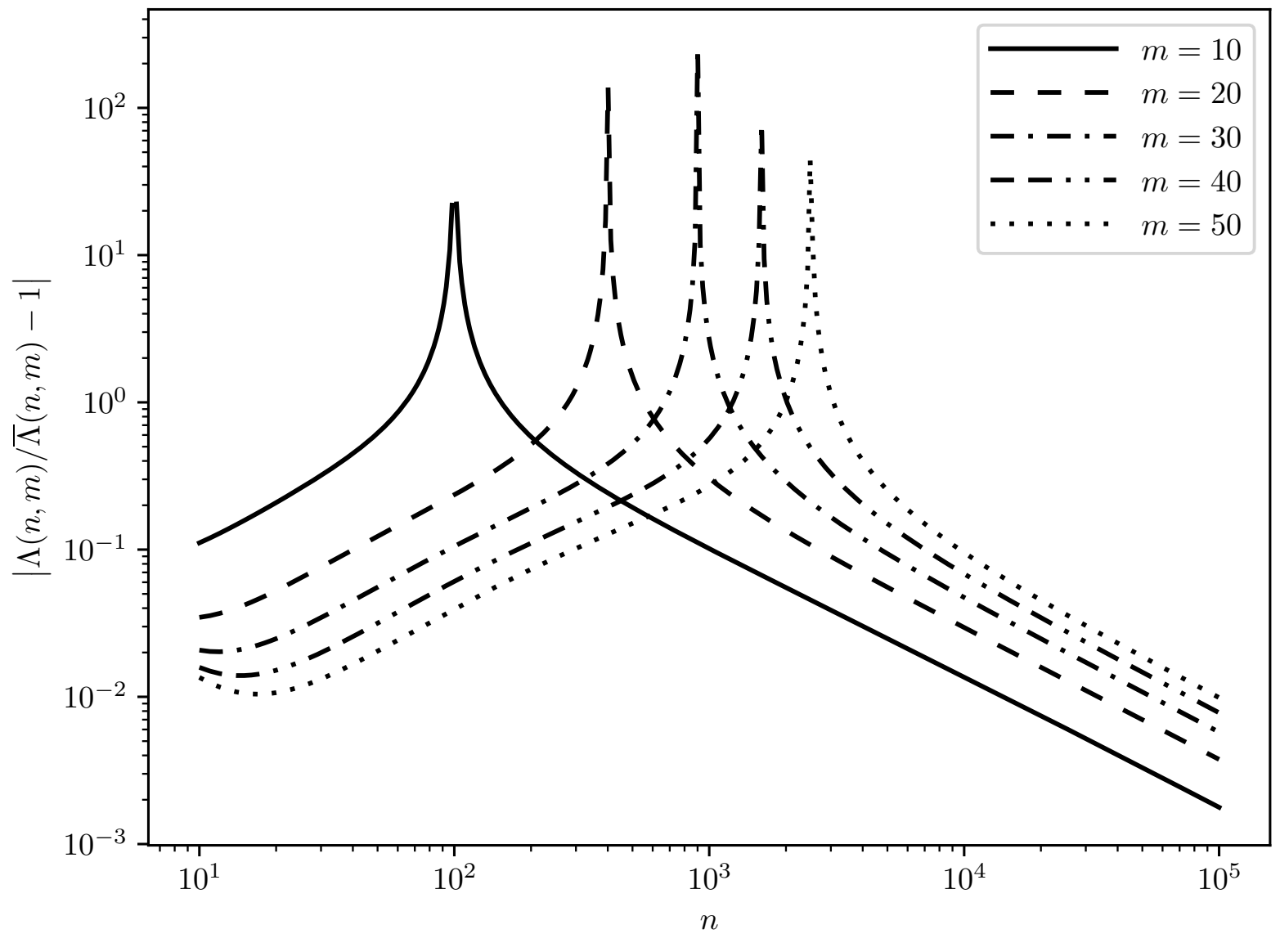

Figure 5 - Comparing predicted Steklov eigenvalues with numerical results. Here ($\Lambda/\overline{\Lambda}-1$) is plotted against $n$ for $m$=10,$\cdots$,50. The spikes occur at roughly $n = m^2$, where the denominator is close to zero.

which may be obtained via numerical quadrature methods. Since the integrand has a singularity at the lower limit of integration $s = 1$, we chose to implement this using tanh-sinh quadrature (also known as "double exponential" quadrature), which is more suitable than other methods (such as Gaussian quadrature) for such problems (c.f. [18], [22] and [23]). After applying numerical root-finding to (4.5), we then simply compute $W = g(H)$.

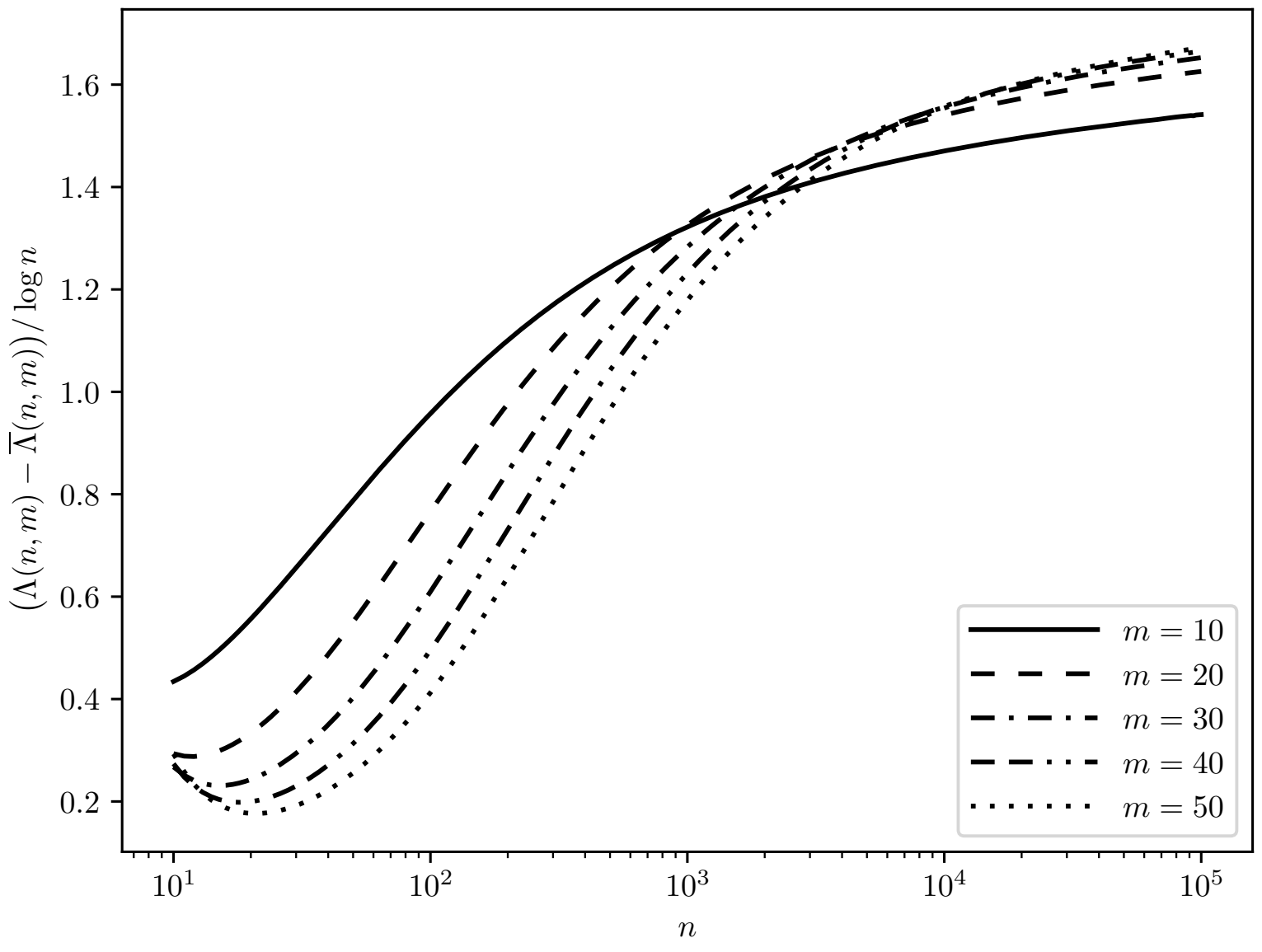

Figure 6 - Comparing predicted Steklov eigenvalues with numerical results. Here $(\Lambda - \overline{\Lambda})/\mathrm{Log}(n)$ is plotted against $n$ for $m$=10,...,50.

Although both methods are equivalent, and both yield identical results, we prefer the latter approach since, not only is numerical quadrature typically more computationally efficient than numerical integration, but tanh-sinh quadrature is also especially amenable to computing with arbitrary-precision arithmetic (c.f. [3]).

Numerical values of $W$ and $H$ are presented in Tables 2 and 3 respectively for $n = 2, \cdots, 100$. We used the `mpmath` [14] arbitrary-precision library in order to validate the results of our (much-faster) floating-point calculations to the level of precision shown here.

**4.3 - Numerical integration for $f$ and $\phi$.** Having obtained numerical values for $W$, we now compute numerical solutions to the initial value problem,

$$f_{xx} = (n-1)(1+f_x^2)f^{-1}, \qquad f(0) = 1, \quad f_x(0) = 0,$$
$$\phi_{xx} = \big(m(m+n-2)(1+f_x^2) - n(n-1)f^2\big)\,\phi, \qquad \phi(0) = 1, \quad \phi_x(0) = 0.$$

over the interval $[0, W]$. As in the first approach presented in Section 4.2, this is carried out using the Dormand and Prince 8(5,3) Runge–Kutta method (specifically, the DOP853 routine described in [13]) with adaptive step-size control. The resulting numerical values of $\phi(W)$ and $\phi_x(W)$ are then used to compute $\Lambda = W\phi_x(W)/\phi(W)$.

Numerical values of $\Lambda$ are given in Table 4 for $m = 2, \cdots, 10$ and $n = 2, \cdots, 20$. We again used the `mpmath` [14] arbitrary-precision library (here, using a Taylor method for arbitrary-precision numerical integration) in order to validate the results of the floating-point calculations to the level of precision displayed in the table. Here the asymptotically linear dependence of $\Lambda$ on $n$ for any fixed $m$ is readily seen by examining the first differences in each of the columns.

| $n$ | $K_0$ | $n$ | $K_0$ | $n$ | $K_0$ | $n$ | $K_0$ | $n$ | $K_0$ |
|---|---|---|---|---|---|---|---|---|---|
| | | 21 | 4 | 41 | 6 | 61 | 7 | 81 | 8 |
| 2 | 2 | 22 | 5 | 42 | 6 | 62 | 7 | 82 | 8 |
| 3 | 2 | 23 | 5 | 43 | 6 | 63 | 7 | 83 | 8 |
| 4 | 2 | 24 | 5 | 44 | 6 | 64 | 7 | 84 | 8 |
| 5 | 3 | 25 | 5 | 45 | 6 | 65 | 7 | 85 | 8 |
| 6 | 3 | 26 | 5 | 46 | 6 | 66 | 7 | 86 | 8 |
| 7 | 3 | 27 | 5 | 47 | 6 | 67 | 7 | 87 | 8 |
| 8 | 3 | 28 | 5 | 48 | 6 | 68 | 7 | 88 | 8 |
| 9 | 3 | 29 | 5 | 49 | 6 | 69 | 7 | 89 | 8 |
| 10 | 3 | 30 | 5 | 50 | 6 | 70 | 8 | 90 | 8 |
| 11 | 3 | 31 | 5 | 51 | 7 | 71 | 8 | 91 | 9 |
| 12 | 4 | 32 | 5 | 52 | 7 | 72 | 8 | 92 | 9 |
| 13 | 4 | 33 | 5 | 53 | 7 | 73 | 8 | 93 | 9 |
| 14 | 4 | 34 | 5 | 54 | 7 | 74 | 8 | 94 | 9 |
| 15 | 4 | 35 | 6 | 55 | 7 | 75 | 8 | 95 | 9 |
| 16 | 4 | 36 | 6 | 56 | 7 | 76 | 8 | 96 | 9 |
| 17 | 4 | 37 | 6 | 57 | 7 | 77 | 8 | 97 | 9 |
| 18 | 4 | 38 | 6 | 58 | 7 | 78 | 8 | 98 | 9 |
| 19 | 4 | 39 | 6 | 59 | 7 | 79 | 8 | 99 | 9 |
| 20 | 4 | 40 | 6 | 60 | 7 | 80 | 8 | 100 | 9 |

Table 7 - Even harmonic modes with Steklov eigenvalues strictly less that 1 as a function of dimension for $n$=2,$\cdots$,100.

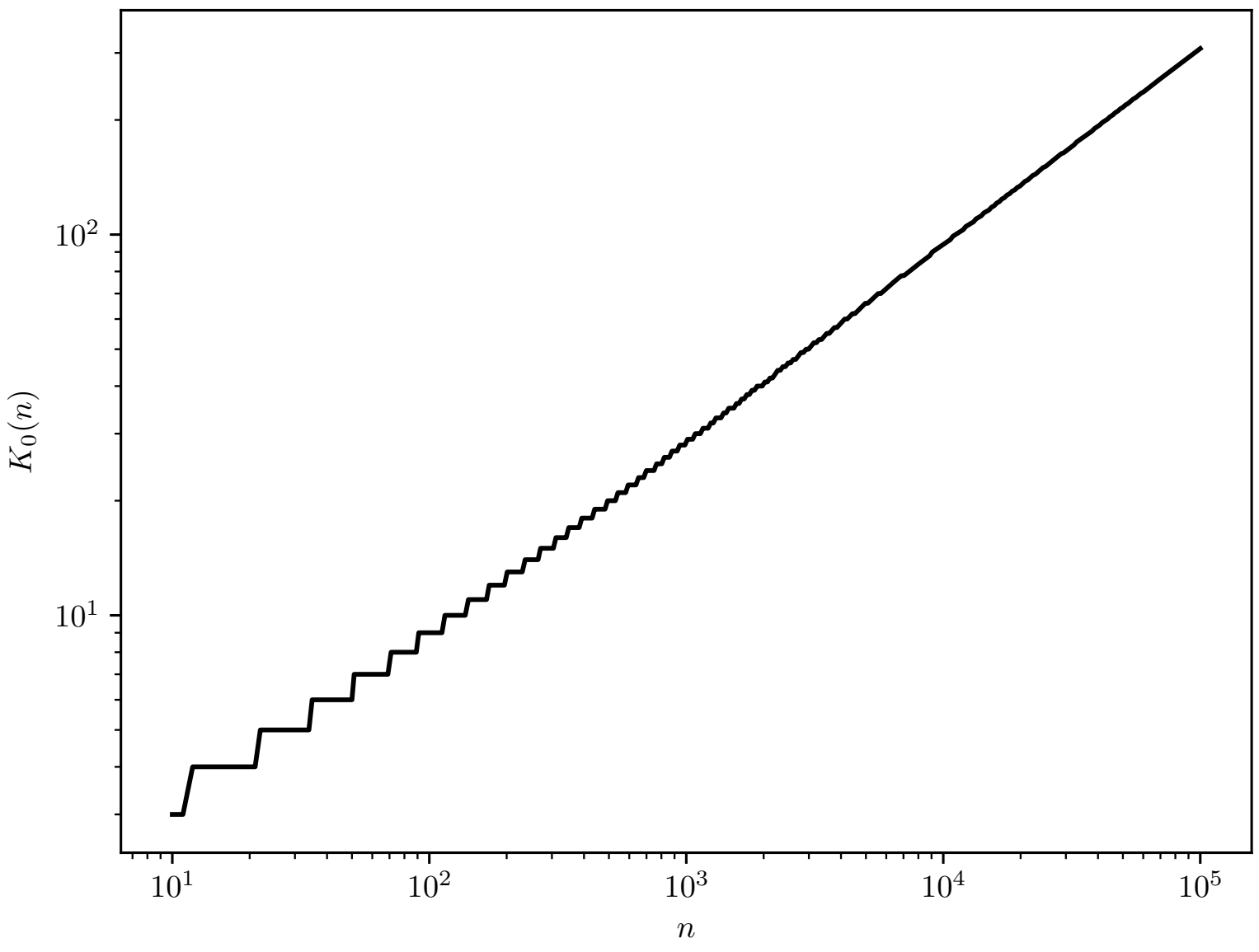

Figure 8 - Even harmonic modes with Steklov eigenvalues strictly less than 1. Here $K_0$ is plotted against $n$ for $n$=10,$\cdots$,$10^5$.

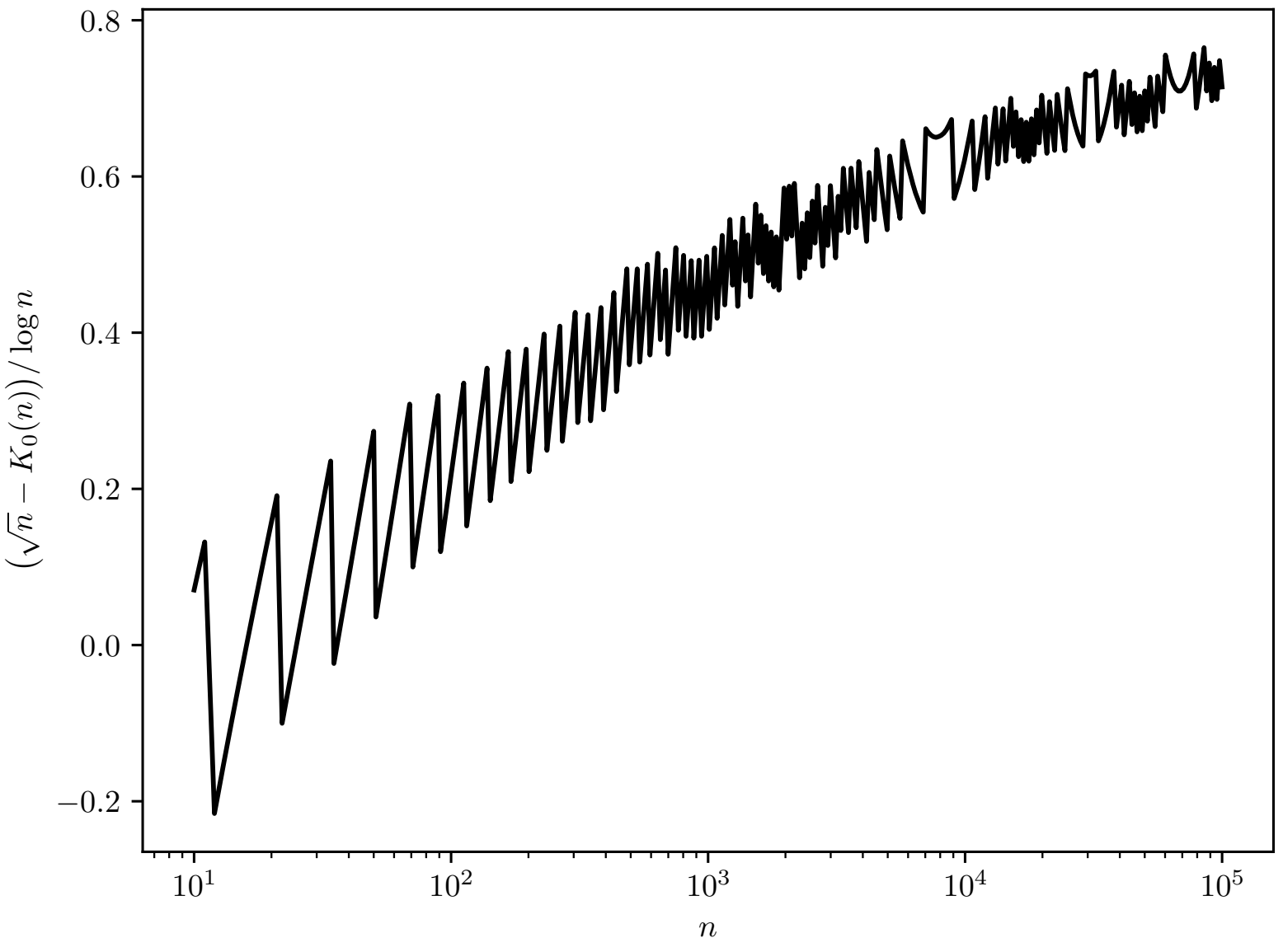

Figure 9 - Comparing predicted values of $\mathrm{K}_0$ with numerical results. Here the rescaled error $(\sqrt{n}-\mathrm{K}_0)/\mathrm{Log}(n)$ is plotted against $n$ for $n$=10,$\cdots$,$10^5$.

We also compare the numerical values of $\Lambda$ with the theoretical estimates of Theorem 3.6.1. Thus, for all $n$ and $m$, we define $\overline{\Lambda} := \overline{\Lambda}_0(n,m)$ by

$$\overline{\Lambda} := m - \frac{n}{m}. \tag{4.6}$$

Figure 5 shows $(\Lambda/\overline{\Lambda} - 1)$ plotted against $n$ for $n = 10, \cdots, 10^5$ and $m = 10, \cdots, 50$. Here we see that the quotient converges to 1 as $n$ tends to infinity, which confirms numerically the highest order term in the asymptotic estimate (3.50). Figure 6 shows $(\Lambda - \overline{\Lambda})/\mathrm{Log}(n)$ plotted against $n$ for $n = 10, \cdots, 10^5$ and $m = 10, \cdots, 50$. Here we see that this quotient remains bounded even for very large values of $n$, which confirms numerically the size of the error in (3.50).

**4.4 - Numerical estimates of the Morse index.** The numerical methods just presented allow us to calculate each Steklov eigenvalue $\Lambda_0(n,m)$ to within an arbitrary error tolerance. By Lemma 2.4.2, we then compute the Morse index $\mathrm{MI}(n)$ by simply counting the number of Steklov eigenvalues strictly less than 1. However, some caution is in order, as the numerical error tolerance must be small enough to resolve the difference between the Steklov eigenvalues and 1. That is, if $k = \mathrm{K}_0(n)$, so that

$$\Lambda_0(n,k-1) < 1 \leq \Lambda_0(n,k),$$

then the numerical error must be less than

$$\mathrm{Min}\left(1 - \Lambda_0(n,k-1), \Lambda_0(n,k) - 1\right),$$

| $n$ | MI | $n$ | MI | $n$ | MI | $n$ | MI | $n$ | MI |
|---|---|---|---|---|---|---|---|---|---|
| | | 21 | 2,003 | 41 | 1,357,511 | 61 | 99,118,657 | 81 | 6,313,511,035 |
| 2 | 4 | 22 | 14,675 | 42 | 1,519,750 | 62 | 108,732,625 | 82 | 6,853,318,716 |
| 3 | 5 | 23 | 17,251 | 43 | 1,697,125 | 63 | 119,110,993 | 83 | 7,432,199,809 |
| 4 | 6 | 24 | 20,151 | 44 | 1,890,670 | 64 | 130,301,601 | 84 | 8,052,482,549 |
| 5 | 21 | 25 | 23,401 | 45 | 2,101,465 | 65 | 142,354,499 | 85 | 8,716,604,821 |
| 6 | 28 | 26 | 27,028 | 46 | 2,330,637 | 66 | 155,322,014 | 86 | 9,427,117,987 |
| 7 | 36 | 27 | 31,060 | 47 | 2,579,361 | 67 | 169,258,818 | 87 | 10,186,690,801 |
| 8 | 45 | 28 | 35,526 | 48 | 2,848,861 | 68 | 184,221,997 | 88 | 10,998,113,413 |
| 9 | 55 | 29 | 40,456 | 49 | 3,140,411 | 69 | 200,271,121 | 89 | 11,864,301,463 |
| 10 | 66 | 30 | 45,881 | 50 | 3,455,336 | 70 | 2,387,548,951 | 90 | 12,788,300,266 |
| 11 | 78 | 31 | 51,833 | 51 | 35,947,198 | 71 | 2,623,427,281 | 91 | 170,212,690,389 |
| 12 | 443 | 32 | 58,345 | 52 | 40,108,069 | 72 | 2,878,995,901 | 92 | 185,035,275,909 |
| 13 | 547 | 33 | 65,451 | 53 | 44,662,465 | 73 | 3,155,605,311 | 93 | 200,974,925,845 |
| 14 | 666 | 34 | 73,186 | 54 | 49,639,591 | 74 | 3,454,679,086 | 94 | 218,103,015,901 |
| 15 | 801 | 35 | 649,573 | 55 | 55,070,247 | 75 | 3,777,716,801 | 95 | 236,494,681,501 |
| 16 | 953 | 36 | 740,260 | 56 | 60,986,885 | 76 | 4,126,297,033 | 96 | 256,228,974,541 |
| 17 | 1,123 | 37 | 840,789 | 57 | 67,423,667 | 77 | 4,502,080,441 | 97 | 277,389,024,991 |
| 18 | 1,312 | 38 | 951,939 | 58 | 74,416,524 | 78 | 4,906,812,925 | 98 | 300,062,207,446 |
| 19 | 1,521 | 39 | 1,074,529 | 59 | 82,003,216 | 79 | 5,342,328,865 | 99 | 324,340,312,726 |
| 20 | 1,751 | 40 | 1,209,419 | 60 | 90,223,393 | 80 | 5,810,554,441 | 100 | 350,319,724,626 |

Table 10 - Morse index as a function of dimension for $n$=2,$\cdots$,100.

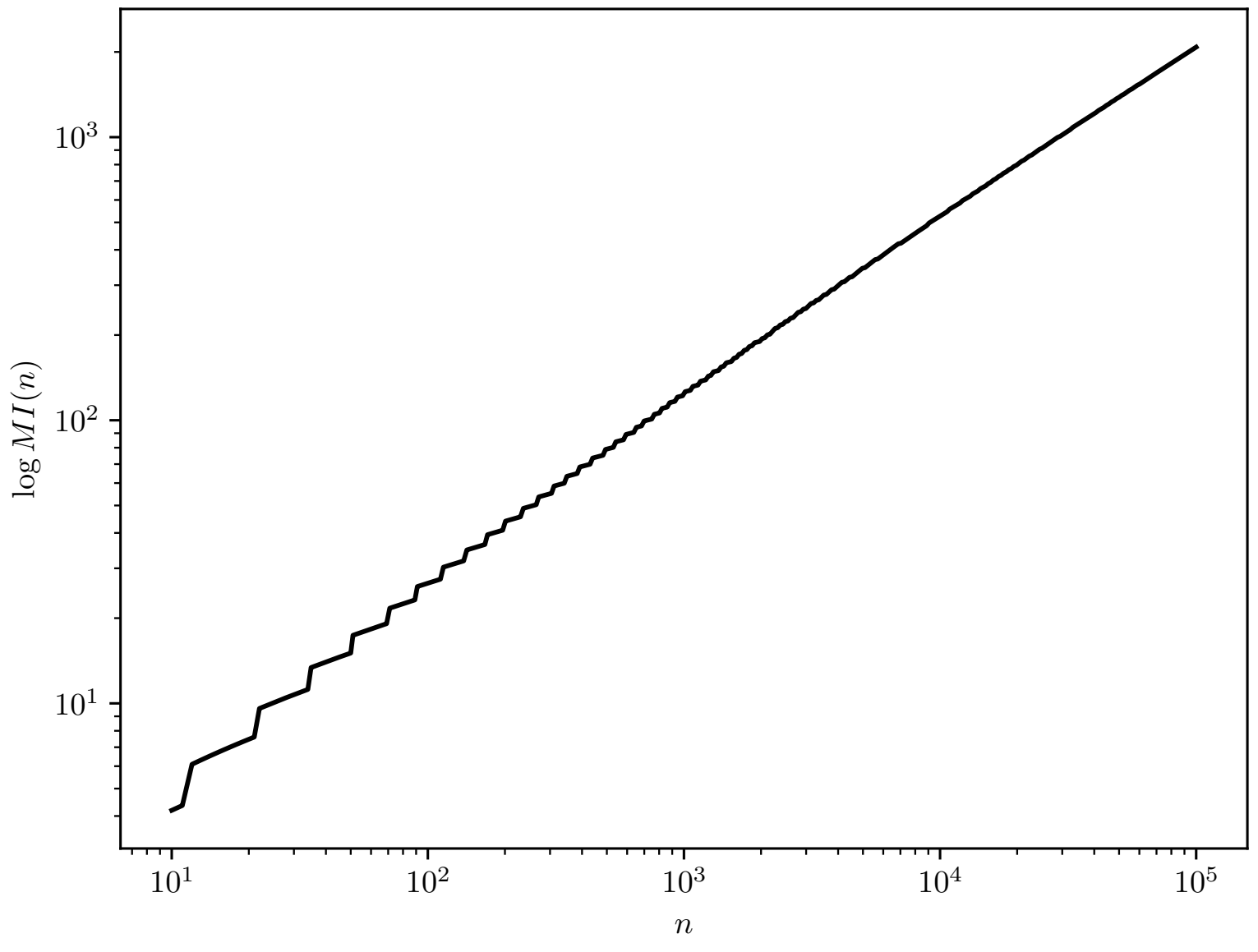

Figure 11 - Morse index as a function of dimension. Here Log(MI) is plotted against $n$ for $n$=10,$\cdots$,$10^5$.

since otherwise an eigenvalue could be less than 1 despite its numerical approximation being greater than 1, or vice versa. However, if the tolerance is indeed taken to be sufficiently small, then we may compute $\mathrm{K}_0(n)$ - and hence $\mathrm{MI}(n)$ - exactly.

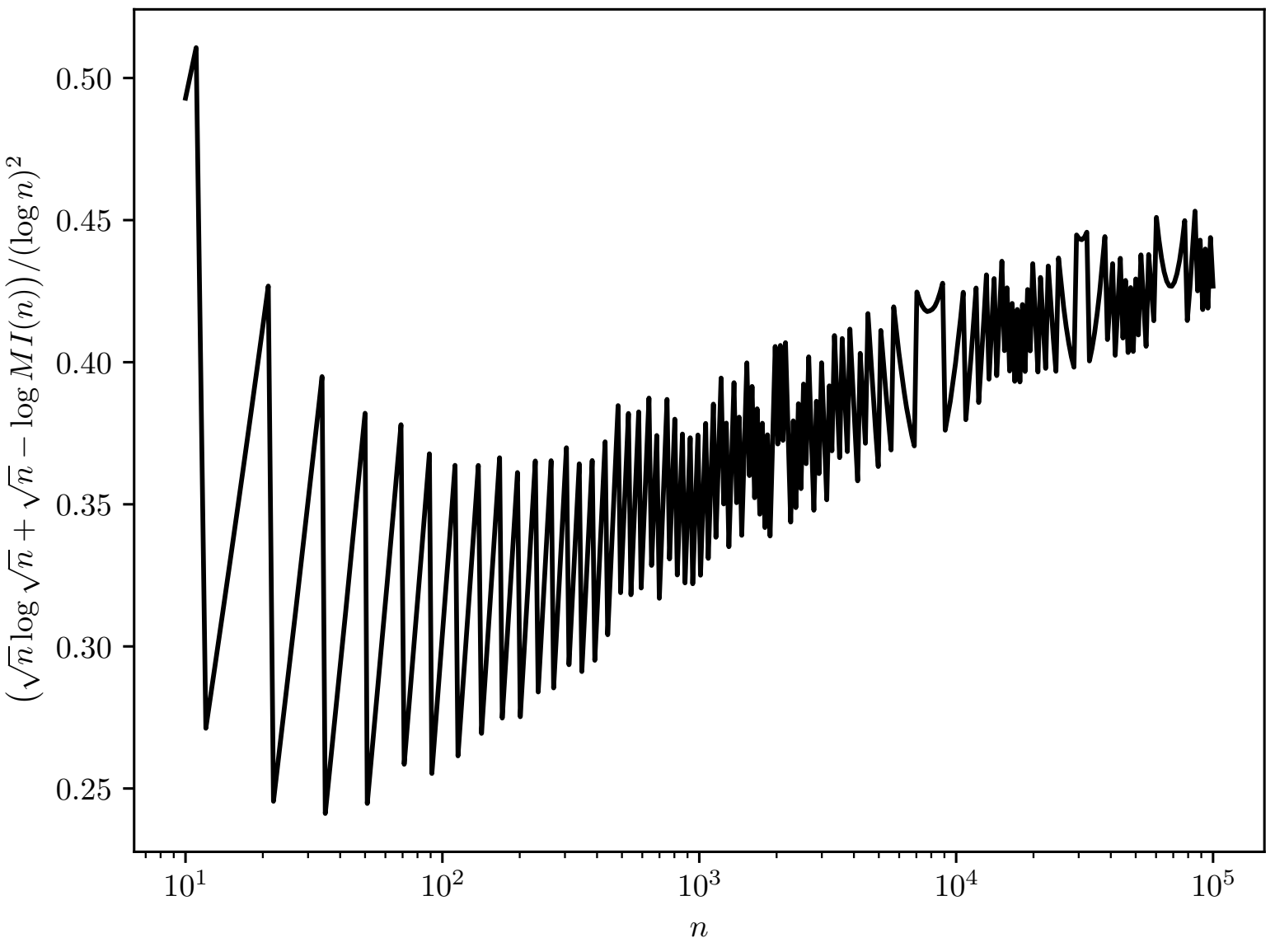

Figure 12 - Comparing predicted values of MI with numerical results. Here the rescaled error $(\sqrt{n}\text{Log}(\sqrt{n})+\sqrt{n}-\text{MI})/\text{Log}(n)^2$ is plotted against $n$ over the range $n$=10,$\cdots$,$10^5$.

Tables 7 and 10 show the values of $\mathrm{K}_0(n)$ and $\mathrm{MI}(n)$, respectively, as computed by this procedure for $n = 2, \cdots, 100$. To ensure that the error tolerance was taken small enough (according to the criteria discussed in the previous paragraph), we computed that the largest Steklov eigenvalue less than 1 was $\Lambda_0(91, 8) \approx 0.99545$, while the smallest Steklov eigenvalue greater than 1 was $\Lambda_0(11, 3) \approx 1.02647$. This was verified using both floating-point arithmetic and arbitrary-precision arithmetic, indicating that the numerical error is much smaller than the tolerance needed to resolve the difference between these eigenvalues and 1.

Denote $K := K_0(n)$ and let $\overline{K} := \overline{K}_0(n, m)$ denote its theoretical estimated value, that is

$$\overline{K} := \overline{K}(n) := \sqrt{n}. \tag{4.7}$$

Figure 8 shows K plotted against $n$ for $n = 10, \cdots, 10^5$. Here we see that the log-log plot is asymptotic to a straight line with gradient 1/2, which confirms numerically the asymptotic estimate (3.52). Figure 9 shows $(K - \overline{K})/\text{Log}(n)$ plotted against $n$ for $n = 10, \cdots, 10^5$. Here we see that the quotient remains bounded even for large values of $n$, which confirms numerically the size of the error in (3.52).

Finally, denote $\mathrm{MI} := \mathrm{MI}(n)$ and let $\overline{\mathrm{MI}}$ denote its theoretical estimated value, that is

$$\overline{\mathrm{MI}} := \overline{\mathrm{MI}}(n) := \text{Exp}\left(\sqrt{n}\text{Log}(\sqrt{n}) + \sqrt{n}\right). \tag{4.8}$$

Figure 11 shows MI plotted against $n$ for $n = 10, \cdots, 10^5$. Likewise, Figure 12 shows $(\text{Log}(\overline{\mathrm{MI}}) - \text{Log}(\mathrm{MI}))/\text{Log}(n^2)$ plotted against $n$ for $n = 10, \cdots, 10^5$. Again, we see that the quotient remains bounded even for large values of $n$, which confirms numerically the size of the error in (3.53).

## A - On the spectra of products of spheres.

For all positive, integer $n$, let $S^n$ denote the unit sphere in $\mathbb{R}^{n+1}$ and, given $1\leq p\leq q$, consider the Cartesian product

$$\Sigma^{p,q}:=\sqrt{\frac{p}{p+q}}S^p\times\sqrt{\frac{q}{p+q}}S^q. \tag{A.1}$$

We readily verify that $\Sigma^{p,q}$ is a minimal hypersurface in $S^{p+q+1}$. The *stability operator* $Q$ of $\Sigma^{p,q}$ is defined to be the second variation of the volume functional with respect to normal perturbations of this hypersurface in $S^{p+q+1}$, and its *Morse index* is defined to be the number of strictly negative eigenvalues of this operator counted with multiplicity.

**Theorem A.0.1**

*If MI$(p,q)$ denotes the Morse index of $\Sigma^{p,q}$ considered as a minimal hypersurface in $S^{p+q+1}$, then*

$$MI(p,q)=p+q+3.$$

**Proof:** Consider the natural parametrisation of $\Sigma^{p,q}$ by $S^p\times S^q$. Recall that the Jacobi operator of $\Sigma^{p,q}$ considered as a minimal hypersurface in $S^{p+q+1}$ is given by

$$\mathrm{J}^{p\cdot q}\phi:=-\Delta\phi-\big((p+q)+\mathrm{Tr}(A^2)\big)\,\phi,$$

where $\Delta^{\Sigma^{p,q}}$ here denotes the Laplace-Beltrami operator of $\Sigma^{p,q}$ and $A$ denotes its shape operator. We readily verify that the principal curvatures of $\Sigma^{p,q}$ in directions tangent to $S^p\times\{0\}$ and $\{0\}\times S^q$ are equal to $\pm\sqrt{q/p}$ and $\mp\sqrt{p/q}$ respectively. From this it follows that

$$\mathrm{J}^{p,q}\phi=-\frac{(p+q)}{p}\Delta^{S^p}\phi-\frac{(p+q)}{q}\Delta^{S^q}\phi-2(p+q)\phi,$$

where, for all $n$, $\Delta^{S^n}$ denotes the Laplace Beltrami operator of $S^n$. As in Section 2.3, for all $(i,j)$, the component of this operator in the $i$'th harmonic mode over $S^p$ and the $j$'th harmonic mode over $S^q$ is then given by

$$\mathrm{J}^{p,q}_{i,j}\phi=\frac{(p+q)}{pq}\left(iq(i+p-1)+jp(j+q-1)\right)\phi-2(p+q)\phi,$$

and we henceforth identify this operator with its unique eigenvalue. Observe that if $i'\geq i$ and if $j'\geq j$, then

$$\mathrm{J}^{p,q}_{i',j'}\geq\mathrm{J}^{p,q}_{i,j},$$

with equality if and only if $i'=i$ and $j'=j$. Furthermore, we readily verify by direct calculation that*

$$\mathrm{J}^{p,q}_{1,1}=0,$$

* In fact, functions in the $(1,1)$ harmonic mode are precisely the Jacobi fields of infinitesimal rotations of the ambient space $S^{p+q+1}$.

and so $\mathrm{J}^{p,q}_{i,j}$ is negative only if one of $i$ or $j$ is equal to 0. However, in this case, we have

$$\mathrm{J}^{p,q}_{i,0} = \frac{(p+q)}{p}\left(i^2 + i(p-1) - 2p\right), \text{ and}$$
$$\mathrm{J}^{p,q}_{0,j} = \frac{(p+q)}{q}\left(j^2 + j(q-1) - 2q\right),$$

so that the only harmonic modes over which J is negative are $(1,0)$, $(0,1)$ and $(0,0)$. Since these modes have respective degeneracies $(p+1)$, $(q+1)$ and 1, the result follows. □

## B - Bibliography.